\documentclass{m2an}
\usepackage{epsfig}

\usepackage{amsfonts}
\usepackage{amsmath}
\usepackage{algorithm}
\usepackage{algorithmic}
\usepackage{multirow}

\graphicspath{{figures/}}

\newcommand{\Xe}{{X_{\rm e}}}
\newcommand{\Ye}{{Y_{\rm e}}}
\newcommand{\Ze}{{Z_{\rm e}}}
\newcommand{\Xpe}{{X'_{\rm e}}}
\newcommand{\Ype}{{Y'_{\rm e}}}

\newcommand{\Ba}{{\rm Ba}}
\newcommand{\Br}{{\rm Br}}
\newcommand{\LB}{{\rm LB}}
\newcommand{\UB}{{\rm UB}}
\newcommand{\cN}{{\mathcal N}}
\newcommand{\cD}{{\mathcal D}}

\newcommand{\bNmax}{{\mathbb{N}_{\rm max}}}

\begin{document}
\title{Reduced Basis {\em A Posteriori} Error Bounds for Symmetric Parametrized Saddle Point Problems}\thanks{This work was supported by the Deutsche Forschungsgemeinschaft (German Research Foundation) through grant GSC 111.}
\author{Anna-Lena Gerner}\address{Aachen Institute for Advanced Study in Computational Engineering Science (AICES), RWTH Aachen University, Schinkelstr.~2, 52062 Aachen, Germany ({\tt gerner@aices.rwth-aachen.de})}
\author{Karen Veroy-Grepl}\address{Aachen Institute for Advanced Study in Computational Engineering Science (AICES) and Faculty of Civil Engineering, RWTH Aachen University, Schinkelstr.~2, 52062 Aachen, Germany ({\tt veroy@aices.rwth-aachen.de})}
%
%
\begin{abstract} 
This paper directly builds upon previous work in \cite{Gerner:2011fk}, where we introduced new reduced basis {\em a posteriori} error bounds for parametrized saddle point problems based on Brezzi's theory. We here sharpen these estimates for the special case of a symmetric problem. Numerical results provide a direct comparison with former approaches and quantify the superiority of the new developed error bounds in practice: Effectivities now decrease significantly; consequently, the proposed methods provide accurate reduced basis approximations at much less computational cost.
\end{abstract}
%
%
\subjclass{65N12, 65N15, 65N30, 76D07}
\keywords{Saddle point problem; Stokes equations; incompressible fluid flow; model order reduction; reduced basis method; {\em a posteriori} error bounds; inf-sup condition}
\maketitle


\section*{Introduction}
The reduced basis (RB) method is a model order reduction approach that permits the efficient yet reliable approximation of input-output relationships induced by parametrized partial differential equations. In contrast to generic discretization techniques where approximation spaces are not correlated to the physical properties of the system, the method recognizes that the solutions to a parametrized partial differential equation are not arbitrary members of the infinite-dimensional solution space but rather reside or evolve on a much lower-dimensional manifold. Exploitation of this low-dimensionality is the key idea of the RB approach.

Designed for the real-time and many-query context of parameter estimation, optimization, and control, the method provides rapidly convergent and computationally efficient approximations equipped with practicable and rigorous error bounds. Built upon a high-fidelity ``truth'' finite element discretization, the RB approximation is defined as a Galerkin projection onto a low-dimensional subspace that focuses on the solution manifold induced by the parametrized partial differential equation. 
The error in the RB approximation is then measured relative to the ``truth'' problem formulation and can be quantified by rigorous and inexpensive {\em a posteriori} error bounds. 
To the method's key features also belongs an Offline-Online computational strategy that enables the highly efficient (Online) computation of both RB approximations and error bounds for any parameter query at the expense of increased pre-processing (Offline) cost. Finally, RB approximations and error bounds are intimately linked through a greedy sampling approach, in which the (Online-)inexpensive error bounds are used to construct the RB approximation spaces more optimally.  

This paper directly builds upon our work in \cite{Gerner:2011fk}, where we introduced new RB {\em a posteriori} error bounds for parametrized saddle point problems based on Brezzi's theory \cite{Brezzi:1974vn,Brezzi:1991fk}. In contrast to former approaches based on Babu\v{s}ka's theory for noncoercive problems \cite{Babuska:1971fk,Manzoni:2012uq,Veroy:2003fk}, 
these do not involve the highly expensive estimation of the Babu\v{s}ka inf-sup stability constants but only much less expensive calculations; as separate upper bounds $\Delta^u$ and $\Delta^p$ for the errors in the RB approximations for the primal variable $u$ and the Lagrange multiplier $p$, respectively, they moreover enable the systematic estimation of engineering outputs depending on either of the two. Numerical results showed that the bounds $\Delta^u$ were reasonably sharp and thus very useful in applications. However, $\Delta^p$ overestimated the actual error in the RB approximations for $p$ rather pessimistically; error estimates $\Delta^\Ba$ based on Babu\v{s}ka's theory here achieved better results in terms of sharpness.

In this paper, we now investigate {\em symmetric} parametrized saddle point problems. We find that in this special case, {\em a posteriori} error bounds $\Delta^u,\Delta^p$ developed in \cite{Gerner:2011fk} may be improved considerably. Considering the Stokes model problem as before, numerical results quantify this behavior in practice through a direct comparison with former techniques: Effectivities for both $\Delta^u$ and $\Delta^p$ now decrease significantly and perform in neither case worse than effectivities associated with $\Delta^\Ba$. 

The paper is organized as follows: In \S \ref{s:general_problem_statement}, we recall the general formulation of a symmetric parametrized saddle point problem together with its ``truth'' approximation. Section \ref{s:RB_method} then describes the RB method with a strong focus on specifics to the symmetric case considered. In \S \ref{ss:Galerkin_projection}, we define the RB approximation as the Galerkin projection onto a low-dimensional RB approximation space. {\em A priori} as well as {\em a posteriori} error estimates shall be specialized to the symmetric case in \S \ref{ss:RB_errbnds}. RB approximations and {\em a posteriori} error bounds can be computed (Online-)efficiently as summarized in \S \ref{ss:offline_online}. This enables us to employ adaptive sampling processes for constructing computationally efficient RB approximation spaces, which shall be outlined in \S \ref{ss:construction_RB_spaces}. In \S \ref{s:numerical_results}, numerical results demonstrate the performance of the improved {\em a posteriori} error bounds in practice. Finally, in \S \ref{s:conclusion}, we give some concluding remarks.

\section{General Problem Statement} \label{s:general_problem_statement}

\subsection{Formulation} \label{ss:pspp}

We shall briefly recall the setting introduced in \cite{Gerner:2011fk}. Let $\Xe$ and $\Ye$ be two Hilbert spaces with inner products $(\cdot,\cdot)_{\Xe}$, $(\cdot,\cdot)_{\Ye}$ and associated norms $\|\cdot\|_{\Xe} = \sqrt{(\cdot,\cdot)_{\Xe}}$, $\|\cdot\|_{\Ye}=\sqrt{(\cdot,\cdot)_{\Ye}}$, respectively.\footnote{Here and in the following, the subscript e denotes ``exact''.} We define the product space $\Ze \equiv \Xe \times \Ye$, with inner product $(\cdot,\cdot)_{\Ze} \equiv (\cdot,\cdot)_{\Xe} + (\cdot,\cdot)_{\Ye}$ and norm $\|\cdot\|_{\Ze} = \sqrt{(\cdot,\cdot)_{\Ze}}$. The associated dual spaces are denoted by $\Xpe$, $\Ype$, and $Z'_{\rm e}$.

Furthermore, let $\cD \subset \mathbb{R}^{n}$ be a prescribed $n$-dimensional, compact parameter set.  For any parameter $\mu \in \cD$, we then consider the continuous bilinear forms $a(\cdot,\cdot;\mu):\Xe\times \Xe\rightarrow \mathbb{R}$ and $b(\cdot,\cdot;\mu):\Xe\times \Ye \rightarrow \mathbb{R}$,\footnote{For clarity of exposition, we suppress the obvious requirement of nonzero elements in the denominators.} 
\begin{align} \label{eq:gamma_a_e}
\gamma_a^{\rm e}(\mu) & \equiv \sup_{u\in \Xe} \sup_{v\in \Xe}\frac{a(u,v;\mu)}{\|u\|_{\Xe}\|v\|_{\Xe}} < \infty \quad\forall\;\mu\in\cD,\\ \label{eq:gamma_b_e}
\gamma_b^{\rm e}(\mu) &\equiv \sup_{q\in \Ye} \sup_{v\in \Xe} 
\frac{b(v,q;\mu)}{\|q\|_{\Ye}\|v\|_{\Xe}} < \infty\quad\forall\;\mu\in\cD.
\end{align}
We moreover assume that $a(\cdot,\cdot;\mu)$ is coercive on $\Xe$,
\begin{equation} \label{eq:alpha_a_e}
\alpha_a^{\rm e}(\mu) \equiv \inf_{v\in \Xe} \frac{a(v,v;\mu)}{\|v\|^2_{\Xe}} > 0 \quad\forall\;\mu\in\cD,
\end{equation}
and that $b(\cdot,\cdot;\mu)$ satisfies the inf-sup condition
\begin{equation} \label{eq:brezzi_infsup_e}
\beta_\Br^{\rm e}(\mu) \equiv \inf_{q\in \Ye} \sup_{v\in \Xe} \frac{b(v,q;\mu)}{\|q\|_{\Ye}\|v\|_{\Xe}}>0\quad\forall\; \mu\in\cD.
\end{equation}
 
We now consider the following variational problem: For any given $\mu\in \cD$, we find $(u_{\rm e}(\mu),p_{\rm e}(\mu)) \in \Xe \times \Ye$ such that
\begin{align} \label{eq:exact} \begin{split}
a(u_{\rm e}(\mu),v;\mu) + b(v,p_{\rm e}(\mu);\mu) &= f(v;\mu) \quad\forall\;v\in \Xe,\\
b(u_{\rm e}(\mu),q;\mu) &= g(q;\mu) \quad\, \forall\;q\in \Ye,
\end{split} 
\end{align}
where $f(\cdot;\mu)$ and $g(\cdot;\mu)$ are bounded linear functionals in $\Xpe$ and $\Ype$, respectively.
From the results of Brezzi  \cite{Brezzi:1974vn,Brezzi:1991fk}, it is well-known that under the assumptions (\ref{eq:gamma_a_e}), (\ref{eq:gamma_b_e}), (\ref{eq:alpha_a_e}), and (\ref{eq:brezzi_infsup_e}), the above problem (\ref{eq:exact}) is well-posed and has a unique solution for any $f(\cdot;\mu)\in \Xpe$, $g(\cdot;\mu)\in \Ype$.

In contrast to this very general setting considered in \cite{Gerner:2011fk}, we here additionally assume that the bilinear form $a(\cdot,\cdot;\mu)$ is {\em symmetric} for any $\mu\in\cD$. As a continuous, symmetric, and coercive bilinear form, $a(\cdot,\cdot;\mu)$ then defines an inner product on $\Xe$ for any $\mu\in\cD$ with an associated norm $\|\cdot\|_{\Xe,\mu}\equiv\sqrt{a(\cdot,\cdot;\mu)}$ equivalent to $\|\cdot\|_\Xe$.
We note (see, e.g., \cite{Brezzi:1991fk,Raviart:1986uq}) that the solution $(u_{\rm e}(\mu), p_{\rm e}(\mu))$ to (\ref{eq:exact}) then corresponds to a saddle point of the Lagrangian functional
\begin{equation*}
\mathcal L(v,q;\mu) \equiv {\textstyle \frac{1}{2}} \, a(v,v;\mu) + b(v,q;\mu) - f(v;\mu) - g(q;\mu), \quad\;(v,q)\in Z_{\rm e}.
\end{equation*}


\subsection{Truth Approximation} \label{ss:truth}

We now introduce a high-fidelity ``truth'' approximation upon which our RB approximation will subsequently be built.  To this end, let $X$ and $Y$ denote finite-dimensional subspaces of $\Xe$ and $\Ye$, respectively.  We define the product space $Z \equiv X \times Y$ and denote by $\cN$ the dimension of $Z$. We emphasize that the dimension $\cN$ is typically very large. These ``truth'' approximation subspaces inherit the inner products and norms of the exact spaces: $(\cdot,\cdot)_X \equiv (\cdot,\cdot)_{\Xe}$, $\|\cdot\|_X\equiv \|\cdot\|_{\Xe}$, $(\cdot,\cdot)_Y \equiv (\cdot,\cdot)_{\Ye}$, $\|\cdot\|_Y\equiv \|\cdot\|_{\Ye}$, and $(\cdot,\cdot)_Z \equiv (\cdot,\cdot)_{\Ze}$, $\|\cdot\|_Z\equiv \|\cdot\|_{\Ze}$.

Clearly, the continuity and coercivity properties (\ref{eq:gamma_a_e}), (\ref{eq:gamma_b_e}), and (\ref{eq:alpha_a_e}) are passed on to the ``truth'' approximation spaces,
\begin{align} \label{eq:gamma_a}
\gamma_a(\mu) &\equiv \sup_{u\in X} \sup_{v\in X} \frac{a(u,v;\mu)}{\|u\|_X \|v\|_X} < \infty \quad\forall\;\mu\in\cD,\\ \label{eq:gamma_b}
\gamma_b(\mu) &\equiv \sup_{q\in Y} \sup_{v\in X} \frac{b(v,q;\mu)}{\|q\|_Y \|v\|_X} < \infty \quad\forall\;\mu\in\cD,\\ \label{eq:alpha_a}
\alpha_a(\mu) &\equiv \inf_{v\in X} \frac{a(v,v;\mu)}{\|v\|^2_X} > 0 \quad\forall\;\mu\in\cD,
\end{align}
and so is the inner product $a(\cdot,\cdot;\mu)$ for any $\mu\in\cD$; thus, $\|\cdot\|_{X,\mu}\equiv\|\cdot\|_{\Xe,\mu}$ defines a norm on $X$ that is equivalent to $\|\cdot\|_X$. We further assume that the approximation spaces $X$ and $Y$ are chosen such that they satisfy the Ladyzhenskaya--Babu\v{s}ka--Brezzi (LBB) inf-sup condition (see, e.g., \cite{Brezzi:1991fk})
\begin{equation} \label{eq:LBB}
\beta_\Br(\mu) \equiv \inf_{q\in Y} \sup_{v\in X} \frac{b(v,q;\mu)}{\|q\|_{Y} \|v\|_X} \geq \beta^0_\Br(\mu) > 0 \quad \forall\; \mu\in\cD,
\end{equation}
where $\beta^0_\Br(\mu)$ is a constant independent of the dimension $\cN$.

We now define our ``truth'' approximations to be the Galerkin projections of $u_{\rm e}(\mu) \in \Xe$ and $p_{\rm e}(\mu) \in \Ye$ onto $X$ and $Y$, respectively: Given any $\mu \in \cD$, we find $(u(\mu), p(\mu)) \in X \times Y$ such that
\begin{align} \label{eq:truth} \begin{split}
a(u(\mu),v;\mu) + b(v, p(\mu);\mu) &= f(v;\mu) \quad\forall\; v\in X,\\
b(u(\mu),q;\mu) &= g(q;\mu) \quad\, \forall\; q\in Y.
\end{split} \end{align}
As for the exact problem in \S \ref{ss:pspp}, it follows from (\ref{eq:gamma_a}), (\ref{eq:gamma_b}), (\ref{eq:alpha_a}), and (\ref{eq:LBB}) that the ``truth'' problem (\ref{eq:truth}) has a unique solution for any $f(\cdot;\mu)\in \Xpe$, $g(\cdot;\mu)\in \Ype$. 
The bilinear forms $a(\cdot,\cdot;\mu)$ and $b(\cdot,\cdot;\mu)$ define bounded linear operators $A(\mu):X \rightarrow X'$, $B(\mu):X\rightarrow Y'$ and its transpose $B(\mu)^t: Y  \rightarrow X'$ by 
\begin{align*} 
\langle A(\mu) \, u, v\rangle & =  a(u,v;\mu)  \hspace{15.3ex}\; \forall\; u, v\in X,\\ 
\langle B(\mu) \, v,q \rangle & =  b(v,q;\mu) = \langle B(\mu)^t q,v\rangle \quad \forall\; v\in X,\; \forall\;q\in Y;
\end{align*}
here, $\langle \cdot,\cdot\rangle$ denotes the respective dual pairing. The ``truth'' system (\ref{eq:truth}) can thus be equivalently written as
\begin{align*}
A(\mu) \,u(\mu) + B(\mu)^t \, p(\mu) & = f(\mu) \quad{\rm in }\; X',\\
B(\mu)\, u(\mu) & = g(\mu) \hspace{0.2ex}\quad{\rm in }\; Y',
\end{align*}
where $f(\mu)\equiv f(\cdot;\mu)|_{X}\in X'$ and $g(\mu)\equiv g(\cdot;\mu)|_Y\in Y'$ for all $\mu\in\cD$.

\section{The Reduced Basis Method} \label{s:RB_method}

We now turn to the RB method. We here focus on the development of {\em a priori} and {\em a posteriori} error estimates, which shall be specialized to the symmetric context. Other parts of the methodology such as the Offline-Online computational strategy as well as the construction of effective RB approximation spaces shall only be briefly recalled as they have already been extensively discussed in \cite{Gerner:2011fk}.

\subsection{Galerkin Projection} \label{ss:Galerkin_projection}

Suppose that we are given a set of nested, low-dimensional RB approximation subspaces $X_N\subset X_{N+1}\subset X$ and $Y_N\subset Y_{N+1}\subset Y$, $N\in \bNmax \equiv\{1,\ldots,N_{\rm max}\}$. We denote by $N_X$ and $N_Y$ the dimensions of $X_N$ and $Y_N$, respectively, and the total dimension of $Z_N\equiv X_N\times Y_N$ by $N_Z \equiv N_X + N_Y$. The subspaces $X_N$, $Y_N$, and $Z_N$ again inherit all inner products and norms of $X$, $Y$, and $Z$, respectively.  The RB approximation is then defined as the Galerkin projection onto these low-dimensional subspaces: For any given $\mu\in\cD$, we find $u_N(\mu)\in X_N$ and $p_N(\mu)\in Y_N$ such that
\begin{align} \label{eq:rb}\begin{split}
a(u_N(\mu),v_N;\mu) + b(v_N,p_N(\mu);\mu) &= f(v_N;\mu) \quad\forall\;v_N\in X_N,\\
b(u_N(\mu),q_N;\mu) & =  g(q_N;\mu) \quad\, \forall\; q_N\in Y_N.
\end{split} \end{align}
Written in operator notation, the discrete RB system reads
\begin{align}\label{eq:rb_op_mom}
A_N(\mu) \,u_N(\mu) + B_N(\mu)^t \, p_N(\mu) & = f_N(\mu) \quad{\rm in }\; X'_N,\\ \label{eq:rb_op_cont}
B_N(\mu)\, u_N(\mu) & = g_N(\mu) \quad{\rm in }\; Y'_N,
\end{align}
where $f_N(\mu)\equiv f(\cdot;\mu)|_{X_N}\in X'_N$, $g_N(\mu)\equiv g(\cdot;\mu)|_{Y_N} \in Y'_N$, and the bounded linear operators $A_N(\mu): X_N\to X'_N$, $B_N(\mu):X_N\to Y'_N$ and its transpose $B_N(\mu)^t:Y_N\to X'_N$ are given by
\begin{align*} 
 \langle A_N(\mu)\, u_N, v_N \rangle &= a(u_N,v_N;\mu) \qquad\qquad\qquad\qquad\quad\, \forall\; u_N, v_N\in X_N,\\ 
 \langle B_N(\mu)\, v_N, q_N \rangle & =  b(v_N,q_N;\mu) = \langle B_N(\mu)^t\, q_N, v_N\rangle \quad\forall\;v_N\in X_N, \;\forall\;q_N\in Y_N.
\end{align*}
We recall (see \cite{Gerner:2011fk}) that the system (\ref{eq:rb}) is well-posed if and only if the RB approximation spaces $X_N,Y_N$ satisfy the inf-sup condition
\begin{equation} \label{eq:rb_infsup}
\beta_N(\mu) \equiv \inf_{q_N\in Y_N} \sup_{v_N\in X_N} \frac{b(v_N,q_N;\mu)}{\|q_N\|_Y\|v_N\|_X} > 0 \quad\forall\;\mu\in\cD;
\end{equation}
in this case, the pair $(X_N, Y_N)$ is called {\em stable}.

\subsection{Reduced Basis Error Estimation} \label{ss:RB_errbnds}

We here consider {\em a priori} as well as {\em a posteriori} estimates for the errors in the RB approximations.

In this section, we assume that the low-dimensional RB spaces $X_N, Y_N$ are constructed such that for any given parameter $\mu\in\cD$, a solution $(u_N(\mu),p_N(\mu))\in X_N\times Y_N$ to (\ref{eq:rb}) exists (see \cite[\S 2.3]{Gerner:2011fk}). We then denote the errors in the RB approximations $u_N(\mu)\in X_N$, $p_N(\mu)\in Y_N$, and $(u_N(\mu),p_N(\mu))\in Z_N$ with respect to the truth approximations by 
\begin{align} \label{eq:def_errors} 
e^u_N(\mu) \equiv u(\mu) - u_N(\mu) \in X,\quad e^p_N(\mu) \equiv p(\mu) - p_N(\mu) \in Y,\quad
e_N(\mu) \equiv (e^u_N(\mu), e^p_N(\mu)) \in Z.
\end{align}

\subsubsection{{\em A Priori} Error Estimates} \label{sss:apriori_errbnds}

Concerning the rate at which the RB approximations $u_N(\mu)$ and $p_N(\mu)$ converge to the truth approximations $u(\mu)$ and $p(\mu)$, respectively, we can derive the following result. 

\begin{prpstn} \label{prpstn:apriori_errbnds}
For any given $\mu\in \cD$ and $N\in\bNmax$, we have
\begin{equation} \label{eq:apriori_errbnd_u}
\|e^u_N(\mu)\|_X \leq 2\sqrt{\frac{\gamma_a(\mu)}{\alpha_a(\mu)}} \inf_{\substack{v_N\in X_N\\B_N(\mu)v_N = g_N(\mu)}} \|u(\mu)-v_N\|_X + \frac{\gamma_b(\mu)}{\alpha_a(\mu)}\inf_{q_N\in Y_N} \|p(\mu)-q_N\|_Y;
\end{equation}
moreover, if the spaces $X_N, Y_N$ are stable (see \S \ref{ss:Galerkin_projection}), we also obtain
\begin{equation} \label{eq:apriori_errbnd_p}
\|e^p_N(\mu)\|_Y \leq \Bigg[1+\frac{\gamma_b(\mu)}{\beta_N(\mu)}\Bigg(1 + \sqrt{\frac{\gamma_a(\mu)}{\alpha_a(\mu)}}\Bigg) \Bigg] \inf_{q_N\in Y_N} \|p(\mu)-q_N\|_Y
 + 2\frac{\gamma_a(\mu)}{\beta_N(\mu)}\inf_{\substack{v_N\in X_N\\B_N(\mu)v_N = g_N(\mu)}} \|u(\mu)-v_N\|_X.
\end{equation}
\end{prpstn}
\begin{proof} We here use techniques very similar to those presented in \cite{Brezzi:1991fk} for finite element approximations. Take any parameter $\mu\in\cD$ and $N\in\bNmax$. Recall that $a(\cdot,\cdot;\mu)$ defines an inner product on $X$ and the associated norm is denoted by $\|\cdot\|_{X,\mu}$.

By the definition of the RB approximation in \S \ref{ss:Galerkin_projection} as the Galerkin projection of $(u(\mu), p(\mu))$ onto $X_N\times Y_N$, the errors $e^u_N(\mu)$ and $e^p_N(\mu)$ satisfy
\begin{equation} \label{eq:Galerkin_proj_mom}
a(e^u_N(\mu),v_N;\mu) + b(v_N,e^p_N(\mu);\mu) = 0 \quad\forall\;v_N\in X_N.
\end{equation}

First, we prove that (\ref{eq:apriori_errbnd_u}) holds true.
For any $v_N\in X_N$ such that $B_N(\mu)\,v_N = g_N(\mu)$ in $Y'_N$, we have $v_N-u_N(\mu)\in \ker(B_N(\mu))$ and
\begin{align*} 
\|v_N-u_N(\mu)\|_{X,\mu} 
& \leq \sup_{w_N\in \ker(B_N(\mu))} \frac{a(v_N-u_N(\mu),w_N;\mu)}{\|w_N\|_{X,\mu}}\\
 &= \sup_{w_N\in\ker(B_N(\mu))} \frac{a(v_N-u(\mu),w_N;\mu)+a(e^u_N(\mu),w_N;\mu)}{\|w_N\|_{X,\mu}} \\
 &= \sup_{w_N\in \ker(B_N(\mu))} \frac{a(v_N-u(\mu),w_N;\mu)-b(w_N,e^p_N(\mu);\mu)}{\|w_N\|_{X,\mu}},
\end{align*}
where the last equality follows from (\ref{eq:Galerkin_proj_mom}).
For $w_N\in\ker(B_N(\mu))$, $b(w_N,p_N(\mu);\mu) = b(w_N,q_N;\mu) = 0$ holds for all $q_N\in Y_N$. Inserting this in the inequality above yields, for any $q_N\in Y_N$,
\begin{align}\nonumber
\|v_N-u_N(\mu)\|_{X,\mu} &\leq \sup_{w_N\in \ker(B_N(\mu))} \frac{a(v_N-u(\mu),w_N;\mu)-b(w_N,p(\mu)-q_N;\mu)}{\|w_N\|_{X,\mu}}\\ \label{eq:apriori_u_1}
&\leq \|v_N-u(\mu)\|_{X,\mu} + \frac{\gamma_b(\mu)}{\sqrt{\alpha_a(\mu)}} \|p(\mu)-q_N\|_Y,
\end{align}
where the latter is obtained from the Cauchy--Schwarz inequality for the inner product $a(\cdot,\cdot;\mu)$, (\ref{eq:gamma_b}), and (\ref{eq:alpha_a}).
Using the triangle inequality and (\ref{eq:apriori_u_1}),
\begin{equation}\label{eq:apriori_u_2}
\|e^u_N(\mu)\|_{X,\mu} \leq \|u(\mu)-v_N\|_{X,\mu} + \|v_N-u_N(\mu)\|_{X,\mu} \leq 2 \|u(\mu)-v_N\|_{X,\mu} + \frac{\gamma_b(\mu)}{\sqrt{\alpha_a(\mu)}}\|p(\mu)-q_N\|_Y,
\end{equation}
the {\em a priori} stability estimate (\ref{eq:apriori_errbnd_u}) then follows from (\ref{eq:apriori_u_2}), (\ref{eq:gamma_a}), and (\ref{eq:alpha_a}).

We now turn to (\ref{eq:apriori_errbnd_p}). Assuming that $X_N$ and $Y_N$ are stable, the inf-sup condition (\ref{eq:rb_infsup}) provides
\begin{equation} \label{eq:apriori_p_1}
\beta_N(\mu)\|q_N-p_N(\mu)\|_Y \leq \sup_{v_N\in X_N} \frac{b(v_N,q_N-p_N(\mu);\mu)}{\|v_N\|_X} \quad\forall\;q_N\in Y_N.
\end{equation}
For any $v_N\in X_N$ and $q_N\in Y_N$, we moreover have
\begin{equation*}
b(v_N,q_N-p_N(\mu);\mu) =  b(v_N, q_N-p(\mu);\mu) + b(v_N,e^p_N(\mu);\mu) = b(v_N,q_N-p(\mu);\mu) - a(e^u_N(\mu),v_N;\mu),
\end{equation*}
where the last equality follows from (\ref{eq:Galerkin_proj_mom}). Applying this to (\ref{eq:apriori_p_1}), together with the Cauchy--Schwarz inequality for $a(\cdot,\cdot;\mu)$ and (\ref{eq:gamma_a}), we obtain
\begin{equation} \label{eq:apriori_p_2}
\|q_N-p_N(\mu)\|_Y \leq \frac{\gamma_b(\mu)}{\beta_N(\mu)}\|p(\mu)-q_N\|_Y + \frac{\sqrt{\gamma_a(\mu)}}{\beta_N(\mu)} \|e^u_N(\mu)\|_{X,\mu} \quad\forall\;q_N\in Y_N;
\end{equation}
the {\em a priori} error estimate (\ref{eq:apriori_errbnd_p}) thus holds again from the triangle inequality, (\ref{eq:apriori_p_2}), (\ref{eq:apriori_u_2}), and (\ref{eq:gamma_a}).
\end{proof}

Note that Proposition~\ref{prpstn:apriori_errbnds} provides a sharper bound for the error in the RB pressure approximations than given in \cite{Brezzi:1991fk} (see also \cite{Gerner:2011fk}) for the more general, nonsymmetric case; in case of the velocity, the bounds generally cannot be related.

\subsubsection{{\em A Posteriori} Error Estimates} \label{sss:aposteriori_errbnds}

We now show that the rigorous and computationally efficient RB {\em a posteriori} error bounds $\Delta^{u}_N(\mu)$, $\Delta^{p}_N(\mu)$, and $\Delta^\Br_N(\mu)$ presented in \cite{Gerner:2011fk} may be further sharpened in the case of a symmetric problem.

To formulate rigorous and inexpensive upper bounds for the respective errors defined in (\ref{eq:def_errors}), we rely on two sets of ingredients. The first set of ingredients consists of computationally (Online-)efficient lower and upper bounds to the truth continuity, coercivity, and inf-sup constants (\ref{eq:gamma_a}), (\ref{eq:alpha_a}), and (\ref{eq:LBB}),
\begin{align} \label{eq:alpha_gamma_LB}
\setlength\arraycolsep{2pt}
\begin{array}{ccccc}
\gamma_a^\LB(\mu) & \leq & \gamma_a(\mu) & \leq & \gamma_a^\UB(\mu),\\[0.5ex]
\alpha_a^\LB(\mu) & \leq &  \alpha_a(\mu) & \leq & \alpha_a^\UB(\mu),\\[0.5ex]
\beta_\Br^\LB(\mu) & \leq & \beta_\Br(\mu) & \leq & \beta_\Br^\UB(\mu),
\end{array}
\quad\forall\;\mu\in\cD.
\end{align}
The second set of ingredients consists of dual norms of the residuals associated with the RB approximation,
\begin{equation} \label{eq:residual_dual_norms}
 \|r^1_N(\cdot;\mu)\|_{X'} = \sup_{v\in X} \frac{r^1_N(v;\mu)}{\|v\|_X},\quad \|r^2_N(\cdot;\mu)\|_{Y'} = \sup_{q\in Y} \frac{r^2_N(q;\mu)}{\|q\|_Y},
\end{equation}
where, for all $\mu\in\cD$, $r^1_N(\cdot;\mu)\in X'$ and $r^2_N(\cdot;\mu)\in Y'$ are defined as
\begin{align} \label{eq:residual_1}
 r^1_N(v;\mu) &\equiv f(v;\mu) - a(u_N(\mu),v;\mu) - b(v,p_N(\mu);\mu) \quad\forall\; v\in X,\\ \label{eq:residual_2}
 r^2_N(q;\mu) &\equiv g(q;\mu) - b(u_N(\mu),q;\mu) \hspace{19.3ex} \forall \; q\in Y.
\end{align}

We may now formulate {\em a posteriori} error bounds for the respective errors (\ref{eq:def_errors}) in the RB approximation.

\begin{prpstn} \label{prpstn:br_errbnds_sym}
For any given $\mu\in\cD$, \mbox{$N\in\bNmax$}, and $\alpha_a^\LB(\mu)$, $\gamma_a^\UB(\mu)$, $\beta_\Br^\LB(\mu)$ satisfying (\ref{eq:alpha_gamma_LB}), we define
\begin{align} \label{eq:br_errbnd_u_sym}
{\Delta}^{u,\rm sym}_{N}(\mu) &\equiv \frac{\|r^1_N(\cdot;\mu)\|_{X'}}{\alpha^\LB_a(\mu)} + \sqrt{\frac{\gamma^\UB_a(\mu)}{\alpha^\LB_a(\mu)}} \frac{\|r^2_N(\cdot;\mu)\|_{Y'}}{\beta^\LB_\Br(\mu)},\\ \label{eq:br_errbnd_p_sym}
{\Delta}^{p,\rm sym}_{N}(\mu) &\equiv \Bigg(1+\sqrt{\frac{\gamma_a^\UB(\mu)}{\alpha_a^\LB(\mu)}}\Bigg)\frac{\|r^1_N(\cdot;\mu)\|_{X'}}{\beta^\LB_\Br(\mu)} + \frac{\gamma_a^\UB(\mu)}{\beta^\LB_\Br(\mu)}\frac{\|r^2_N(\cdot;\mu)\|_{Y'}}{\beta^\LB_\Br(\mu)}.
\end{align}
Then, $ \Delta^{u,\rm sym}_N(\mu)$ and $ \Delta^{p,\rm sym}_N(\mu)$ are upper bounds for the errors $e^u_N(\mu)$ and $e^p_N(\mu)$ such that
\begin{equation} \label{eq:br_u_p_rig_sym}
\|e^u_N(\mu)\|_X \leq \Delta^{u,\rm sym}_N(\mu)<\Delta^{u}_N(\mu),\quad \|e^p_N(\mu)\|_Y \leq \Delta^{p,\rm sym}_N(\mu)<\Delta^{p}_N(\mu)
\end{equation}
for all $\mu \in\cD$ and $N\in\bNmax$, where $\Delta^u_N(\mu)$ and $\Delta^p_N(\mu)$ are defined as in \cite[Corollary~2.2]{Gerner:2011fk}.
\end{prpstn}
\begin{proof} Let $\mu$ be an arbitrary but fixed parameter in $\cD$ and $N\in\bNmax$. We here proceed as in the proof of \cite[Corollary~2.2]{Gerner:2011fk}, only that we may now exploit the fact that $a(\cdot,\cdot;\mu)$ defines an inner product on $X$; again, recall that the associated norm is denoted by $\|\cdot\|_{X,\mu}$.

By (\ref{eq:residual_1}), (\ref{eq:residual_2}), and (\ref{eq:truth}), the errors $e^u_N(\mu)\in X$ and $e^p_N(\mu)\in Y$ satisfy the equations
 \begin{align} \label{eq:err_mom}
a(e^u_N(\mu),v;\mu) + b(v,e^p_N(\mu);\mu) & = r^1_N(v;\mu) \quad\forall\; v\in X,\\ \label{eq:err_cont}
b(e^u_N(\mu),q;\mu) & =  r^2_N(q;\mu) \quad \forall \; q\in Y.
\end{align}

The error $e^u_N(\mu)\in X$ in the approximation of the primal variable may now be uniquely decomposed into 
$e^u_N(\mu) = \tilde e^0_N(\mu) + \tilde e^\perp_N(\mu)$ where $\tilde e^0_N(\mu)\in \ker(B(\mu))$ and $\tilde e^\perp_N(\mu)\in X$ such that 
\begin{equation}\label{eq:e_perp_orth}
a(\tilde e^\perp_N(\mu),v_0;\mu)=0 \quad\forall\;v_0\in \ker(B(\mu)).
\end{equation}
From (\ref{eq:err_mom}) and (\ref{eq:e_perp_orth}), $\tilde e^0_N(\mu)\in\ker(B(\mu))$ then solves
\begin{equation} \nonumber
a(\tilde e^0_N(\mu),v_0;\mu) = r^1_N(v_0;\mu) - a(\tilde e^\perp_N(\mu),v_0;\mu) = r^1_N(v_0;\mu) \quad\forall\; v_0\in \ker(B(\mu)).
\end{equation}
Setting here $v_0=\tilde e^0_N(\mu)$, we have
\begin{align*}
\|\tilde e^0_N(\mu)\|_{X,\mu}^2 & = r^1_N(\tilde e^0_N(\mu);\mu) \leq \|\tilde e^0_N(\mu)\|_X \sup_{v_0\in\ker(B(\mu))} \frac{r^1_N(v_0;\mu)}{\|v_0\|_X}\\
&\leq \frac{1}{\sqrt{\alpha_a(\mu)}} \|\tilde e^0_N(\mu)\|_{X,\mu} \sup_{v_0\in\ker(B(\mu))} \frac{r^1_N(v_0;\mu)}{\|v_0\|_X},
\end{align*}
where the last inequality follows from (\ref{eq:alpha_a}). Hence, $\tilde e^0_N(\mu)$ is bounded by
\begin{equation}\label{eq:e_0_bnd_sym}
\|\tilde e^0_N(\mu)\|_{X,\mu} \leq \frac{1}{\sqrt{\alpha_a(\mu)}} \sup_{v_0\in\ker(B(\mu))} \frac{r^1_N(v_0;\mu)}{\|v_0\|_X}
\leq \frac{\|r^1_N(\cdot;\mu)\|_{X'}}{\sqrt{\alpha_a(\mu)}}.
\end{equation}
To obtain an upper bound for $\tilde e^\perp_N(\mu)$, we here consider the inf-sup constant
\begin{equation}\label{eq:LBB_tilde}
\tilde\beta(\mu) \equiv \inf_{q\in Y}\sup_{v\in X} \frac{b(v,q;\mu)}{\|q\|_Y \|v\|_{X,\mu}}.
\end{equation}
From (\ref{eq:gamma_a}) and (\ref{eq:alpha_a}), we have
\begin{equation}\label{eq:tilde_beta_equiv}
\frac{\beta_\Br(\mu)}{\sqrt{\gamma_a(\mu)}} \leq\tilde\beta(\mu)\leq \frac{\beta_\Br(\mu)}{\sqrt{\alpha_a(\mu)}},
\end{equation}
and thus particularly $\tilde\beta(\mu)>0$ by the LBB inf-sup condition (\ref{eq:LBB}). Consequently, 
\begin{equation*}
\tilde\beta(\mu) \|v\|_{X,\mu} \leq \sup_{q\in Y}\frac{b(v,q;\mu)}{\|q\|_Y}
\end{equation*}
holds true for any $v\in X$ such that $a(v,v_0;\mu) = 0$ for all $v_0\in\ker(B(\mu))$ (see, e.g., \cite[\S II.1, Proposition~1.2]{Brezzi:1991fk}). Applied to $\tilde e^\perp_N(\mu)$ satisfying (\ref{eq:e_perp_orth}), this yields 
\begin{equation} \label{eq:e_perp_bnd_sym}
\|\tilde e^\perp_N(\mu)\|_{X,\mu}\leq \frac{1}{\tilde\beta(\mu)} \sup_{q\in Y} \frac{b(\tilde e^\perp_N(\mu),q;\mu)}{\|q\|_Y} = \frac{\|r^2_N(\cdot;\mu)\|_{Y'}}{\tilde\beta(\mu)},
\end{equation}
where the equality follows from (\ref{eq:err_cont}) as $B(\mu)\,\tilde e^\perp_N(\mu) = B(\mu)\,e^u_N(\mu)$ in $Y'$. Now, using the triangle inequality, we may finally derive that
\begin{equation} \label{eq:e_u_energy_bnd} 
\|e^u_N(\mu)\|_{X,\mu} \leq \|\tilde e^0_N(\mu)\|_{X,\mu} + \|\tilde e^\perp_N(\mu)\|_{X,\mu} \leq \frac{\|r^1_N(\cdot;\mu)\|_{X'}}{\sqrt{\alpha_a(\mu)}} + \frac{\sqrt{\gamma_a(\mu)}}{\beta_\Br(\mu)}\|r^2_N(\cdot;\mu)\|_{Y'},
\end{equation}
by combining (\ref{eq:e_0_bnd_sym}), (\ref{eq:e_perp_bnd_sym}), and (\ref{eq:tilde_beta_equiv}); the bound (\ref{eq:br_errbnd_u_sym}) thus follows from (\ref{eq:e_u_energy_bnd}), (\ref{eq:alpha_a}), and (\ref{eq:alpha_gamma_LB}).

For the error $e^p_N(\mu)$ in the approximation of the Lagrange multiplier, we obtain from (\ref{eq:LBB}) and (\ref{eq:err_mom}) that
\begin{align} \nonumber
\|e^p_N(\mu)\|_Y &\leq \frac{1}{\beta_\Br(\mu)} \sup_{v\in X} \frac{b(v,e^p_N(\mu);\mu)}{\|v\|_X}
= \frac{1}{\beta_\Br(\mu)} \sup_{v\in X} \frac{r^1_N(v;\mu) - a(e^u_N(\mu),v;\mu)}{\|v\|_X}\\ \label{eq:ep_bounded_by_eu_sym}
&\leq \frac{1}{\beta_\Br(\mu)} \Big(\|r^1_N(\cdot;\mu)\|_{X'}  + \sqrt{\gamma_a(\mu)}  \|e^u_N(\mu)\|_{X,\mu}\Big),
\end{align}
where the last inequality holds by the Cauchy--Schwarz inequality for the inner product $a(\cdot,\cdot;\mu)$ and (\ref{eq:gamma_a}). Together with (\ref{eq:e_u_energy_bnd}) and again (\ref{eq:alpha_gamma_LB}), this leads to (\ref{eq:br_errbnd_p_sym}).

As it is clearly $\alpha_a(\mu)\leq \gamma_a(\mu)$ from the definitions in (\ref{eq:gamma_a}) and (\ref{eq:alpha_a}), we have $\Delta^{u,\rm sym}_N(\mu)<\Delta^{u}_N(\mu)$ and $\Delta^{p,\rm sym}_N(\mu)<\Delta^{p}_N(\mu)$ (see formulations of $\Delta^{u}_N(\mu)$, $\Delta^{p}_N(\mu)$ in \cite[Corollary~2.2]{Gerner:2011fk}); this eventually yields~(\ref{eq:br_u_p_rig_sym}).
\end{proof}


Now, we clearly obtain
\begin{equation}\label{eq:br_errbnd_sym}
 \Delta^{\rm sym}_N(\mu) \equiv \sqrt{\big(\Delta^{u,\rm sym}_N(\mu)\big)^2 + \big(\Delta^{p,\rm sym}_N(\mu)\big)^2}, \quad\mu\in\cD,\;N\in\bNmax,
\end{equation} 
as a rigorous upper bound for the combined error $e_N(\mu)$ such that
\begin{equation*}
\|e_N(\mu)\|_Z \leq \Delta^{\rm sym}_N(\mu)<\Delta^\Br_N(\mu)\quad\forall\;\mu\in\cD,\;\forall\;N\in\bNmax,
\end{equation*}
where $\Delta^\Br_N(\mu)$ is given as in \cite[(2.20)]{Gerner:2011fk}.
Analogous to \cite[Remark~2.3]{Gerner:2011fk}, the associated effectivities are then bounded by
\begin{equation} \label{eq:br_effectivity_sym}
1\leq \eta^{\rm sym}_N(\mu)\equiv \frac{\Delta^{\rm sym}_N(\mu)}{\|e_N(\mu)\|_Z}\leq C(\mu)\big(\gamma_a(\mu)+\gamma_b(\mu)\big) \quad\forall\;\mu\in\cD,\;\forall\;N\in\bNmax,
\end{equation}
where $C(\mu)>0$ is a constant depending on $\alpha^\LB_a(\mu)$, $\gamma_a^\UB(\mu)$, and $\beta^\LB_\Br(\mu)$.


Furthermore, the error $e^u_N(\mu)$ may particularly also be measured in the energy norm $\|\cdot \|_{X,\mu}= \sqrt{a(\cdot,\cdot;\mu)}$. For any given $\mu\in\cD$, $N\in\bNmax$, and $\alpha_a^\LB(\mu)$, $\gamma_a^\UB(\mu)$, $\beta_\Br^\LB(\mu)$ satisfying (\ref{eq:alpha_gamma_LB}), we define
\begin{equation}\label{eq:br_errbnd_u_sym_energy}
\tilde\Delta^{u,\rm sym}_{N}(\mu) \equiv \frac{\|r^1_N(\cdot;\mu)\|_{X'}}{\sqrt{\alpha_a^\LB(\mu)}} + \frac{\sqrt{\gamma_a^\UB(\mu)}}{\beta_\Br^\LB(\mu)} \|r^2_N(\cdot;\mu)\|_{Y'}.
\end{equation}
Then, from (\ref{eq:e_u_energy_bnd}) and (\ref{eq:alpha_gamma_LB}), $\tilde\Delta^{u,\rm sym}_N(\mu)$ is an upper bound for $e^u_N(\mu)$ such that
\begin{equation}\label{eq:br_errbnd_u_sym_energy_rig}
\|e^u_N(\mu)\|_{X,\mu} \leq \tilde\Delta^{u,\rm sym}_N(\mu) \quad\forall\;\mu\in\cD,\;\forall\;N\in\bNmax.
\end{equation}
It directly follows from the proof of Proposition~\ref{prpstn:br_errbnds_sym} that the effectivities $\tilde\eta^{u,\rm sym}_N(\mu)\equiv \tilde\Delta^{u,\rm sym}_N(\mu)/\|e^u_N(\mu)\|_{X,\mu}$ and $\eta^{u,\rm sym}_N(\mu)\equiv \Delta^{u,\rm sym}_N(\mu)/\|e^u_N(\mu)\|_X$ satisfy 
\begin{equation}\label{eq:eff_sym}
\tilde\eta^{u,\rm sym}_N(\mu)\leq \eta^{u,\rm sym}_N(\mu)\quad \forall\; \mu\in\cD, \;\forall\;N\in\bNmax.
\end{equation}

\subsection{Offline-Online Computational Procedure} \label{ss:offline_online}

The efficiency of the RB method relies on an Offline-Online computational decomposition strategy. As this is by now standard, we shall only provide a brief summary at this point and refer the reader to, e.g., \cite{Rozza:2008fv} for further details. The procedure requires that all involved operators can be affinely expanded with respect to the parameter $\mu$. All $\mu$-independent quantities are then formed and stored within a computationally expensive Offline stage, which is performed only once and whose cost depends on the large truth dimension $\cN$. For any given $\mu\in\cD$, the RB approximation is then computed within a highly efficient Online stage; the cost does not depend on $\cN$ but only on the considerably smaller dimension of the RB approximation space. 

The computation of the {\em a posteriori} error bounds clearly consists of two components: the calculation of the residual dual norms (\ref{eq:residual_dual_norms}) and the lower and upper bounds (\ref{eq:alpha_gamma_LB}) to the involved coercivity, continuity, and inf-sup stability constants. The former is again an application of standard RB techniques that can be found in \cite{Rozza:2008fv}. The latter is achieved by a successive constraint method (SCM) proposed in \cite{DBP-Huynh:2007ez}; we refer the reader to \cite{Gerner:2012ag} for details in our saddle point context.

\subsection{Construction of Reduced Basis Approximation Spaces} \label{ss:construction_RB_spaces}

The low-dimensional RB approximation spaces $X_N,Y_N$, $N\in\bNmax$, are constructed by exploiting the parametric structure of the problem: According to the so-called Lagrange approach, basis functions are given by truth solutions associated with several chosen parameter snapshots.  In the case of saddle point problems, additional care must be taken in the construction of stable RB approximation spaces (see \S \ref{ss:Galerkin_projection}). This has been extensively discussed in \cite{Gerner:2011fk}: Stability is achieved through enriching the RB space $X_N$ for the primal variable appropriately. Different strategies may be applied: We may add {\em supremizer} functions \cite{Gerner:2011fk,Rovas:2003ly,Rozza:2007fk} or additional truth solutions \cite{Gerner:2011fk} favoring either the approximations for the primal or the Lagrange multiplier variables, respectively.

Keeping computational cost to a minimum, we aim to construct stable approximation spaces $X_N,Y_N$ that appropriately represent the submanifold associated with the parametric dependence with as few basis functions as possible. To this end, we invoke a greedy sampling process \cite{dahmen10:_greedy,Buffa:kx} in which our rigorous and computationally inexpensive RB {\em a posteriori} error bounds are used to identify truth solutions that are not yet well approximated. In \cite{Gerner:2011fk}, this procedure has been extended to the special needs of saddle point problems. 
Correspondingly, we here again consider the following three sampling processes: the standard greedy sampling process summarized in \cite[Algorithm~1]{Gerner:2011fk} as well as the modified procedures \cite[Algorithm~2]{Gerner:2011fk} and \cite[Algorithm~3]{Gerner:2011fk} where the need for stabilization is recognized adaptively.

\section{Numerical Results} \label{s:numerical_results}

\begin{figure}[!bp]
\epsfig{file=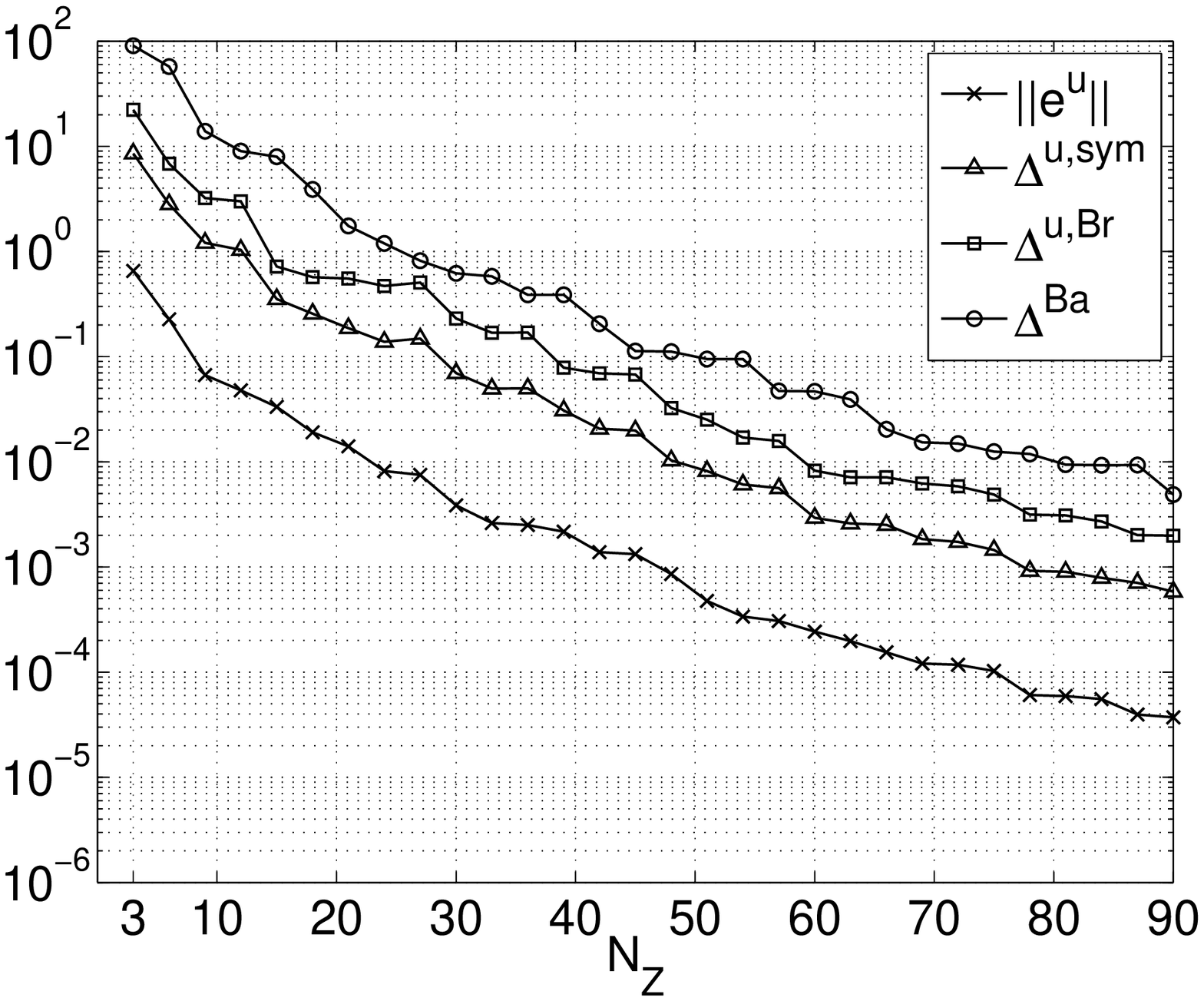,scale=0.3}\hfill
\epsfig{file=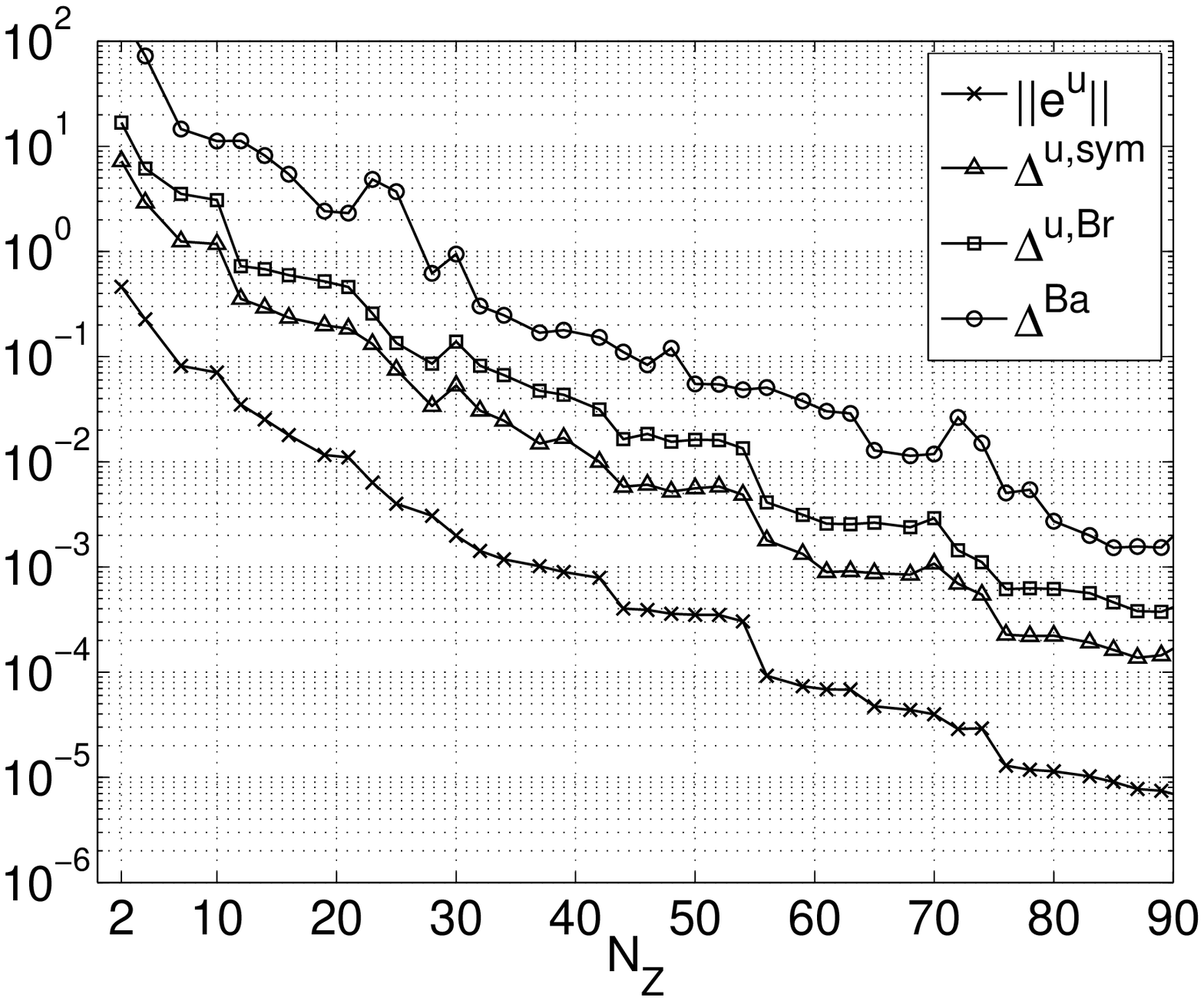,scale=0.3}\hfill
\epsfig{file=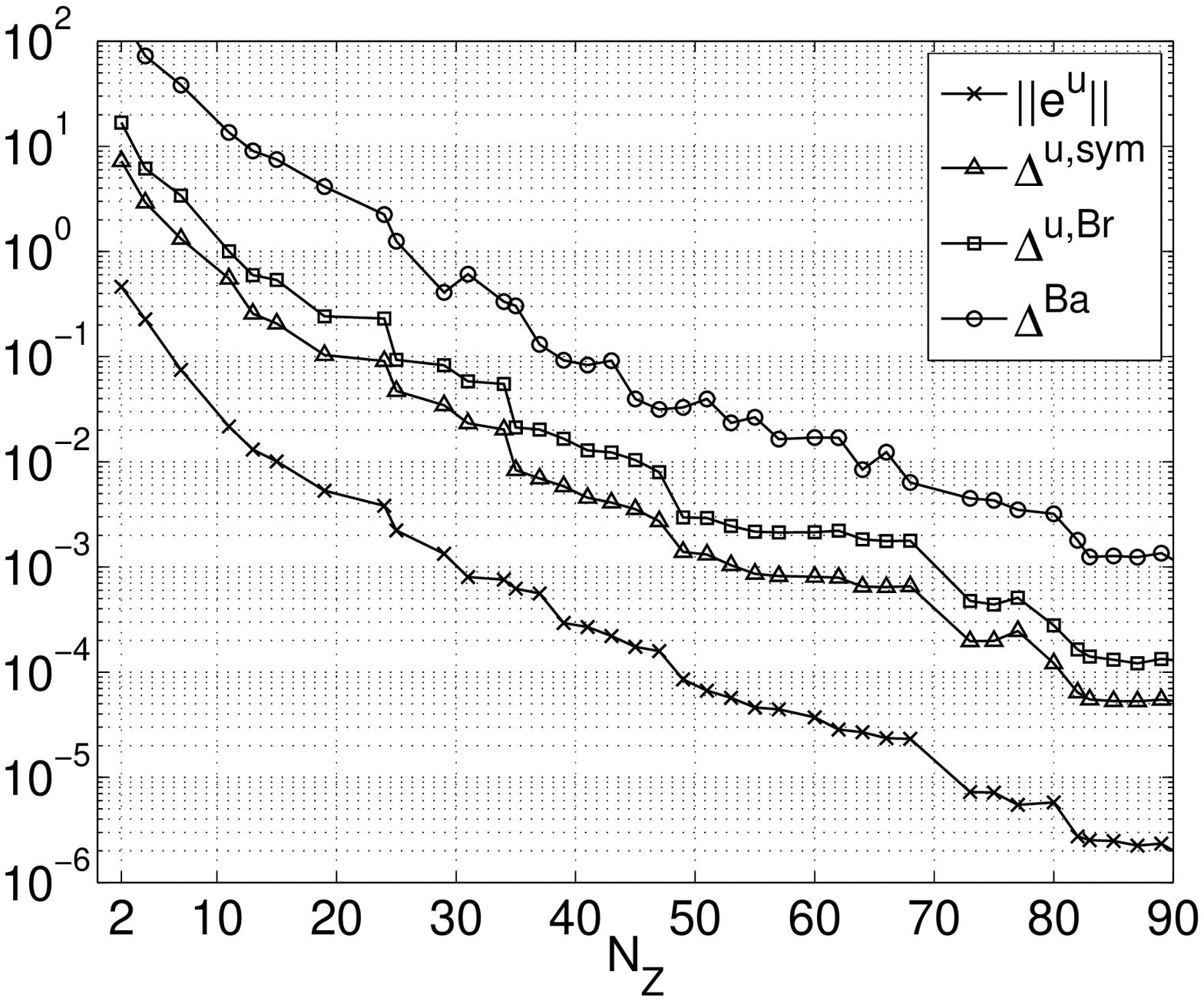,scale=0.3}\hfill
\centerline{(a) \hspace{35ex} (b) \hspace{35ex} (c)}
\caption{Maximum relative error $\|e^u_N(\mu)\|_X/\|u(\mu)\|_X$ and maximum relative error bounds $\Delta^{u,\rm sym}_N(\mu)/\|u(\mu)\|_X$, $\Delta^{u}_N(\mu)/\|u(\mu)\|_X$, and $\Delta^\Ba_N(\mu)/\|u(\mu)\|_X$ shown as functions of $N_Z$ for (a)~Algorithm~1, (b)~Algorithm~2, and (c)~Algorithm~3 (see \S \ref{ss:construction_RB_spaces}); the maximum is taken over 25 parameter values; the computation of the error bounds is based on the exact constants involved.}
\label{fig:u_error}
\end{figure}

\begin{figure}[htp]
\epsfig{file=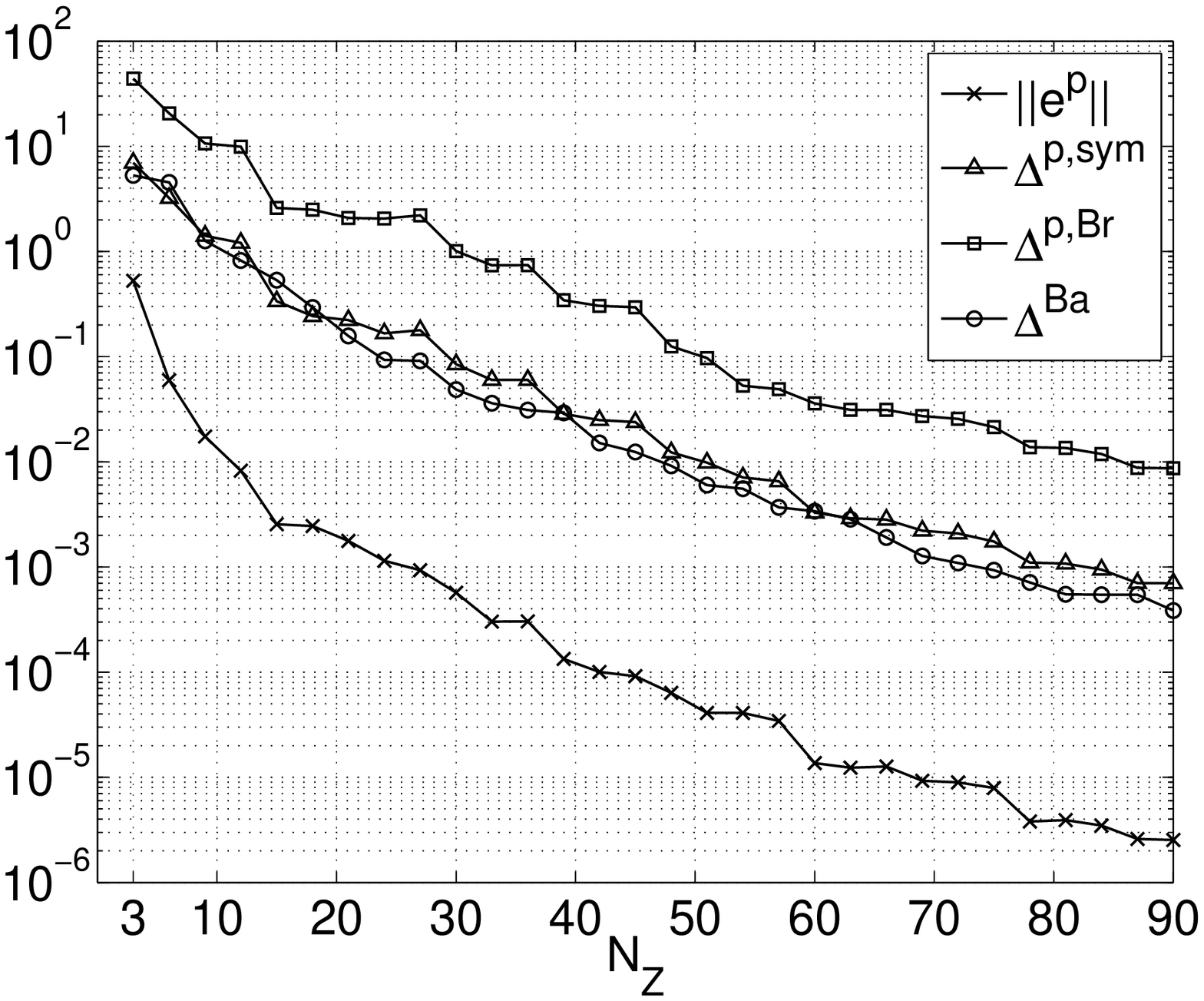,scale=0.3}\hfill
\epsfig{file=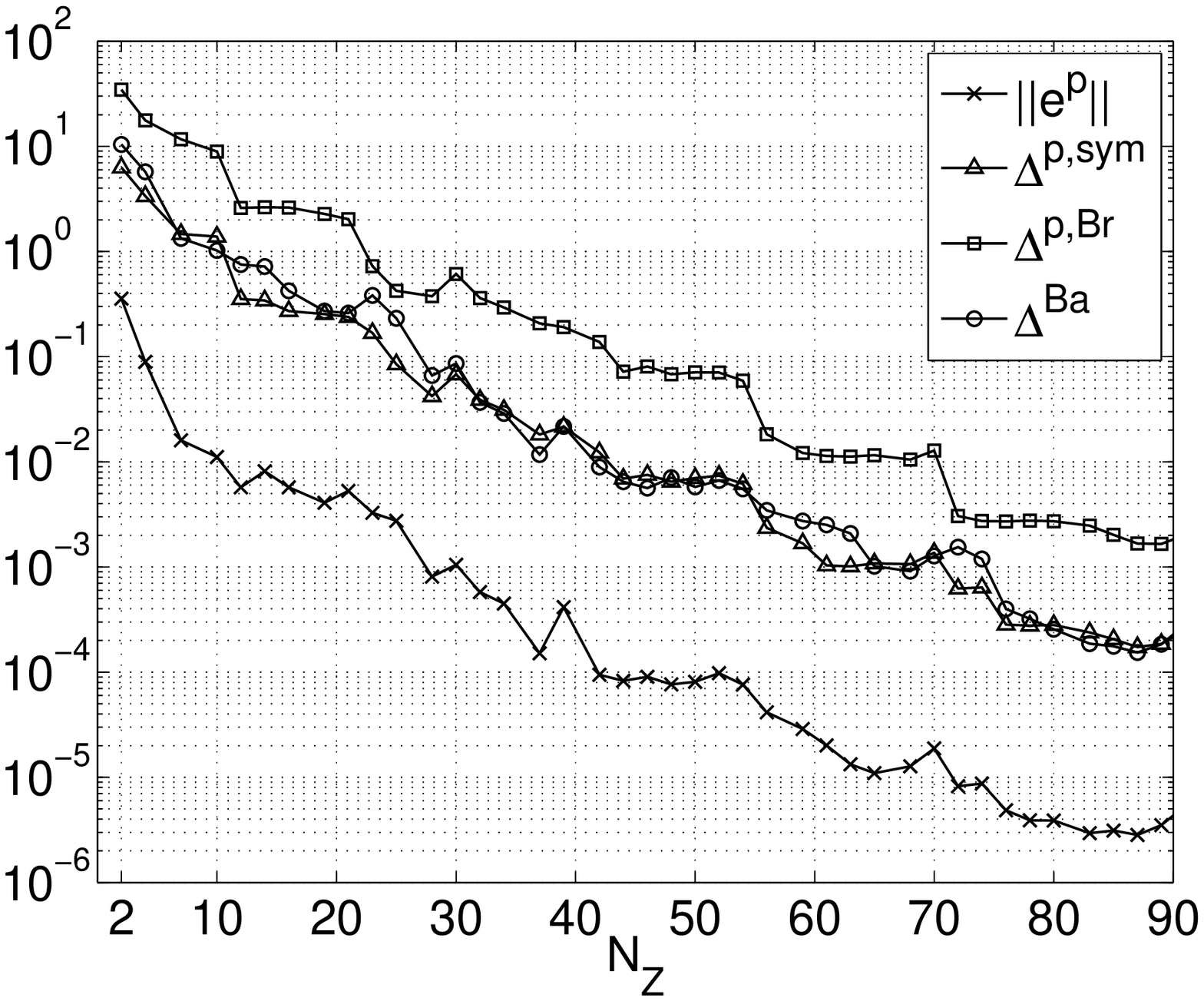,scale=0.3}\hfill
\epsfig{file=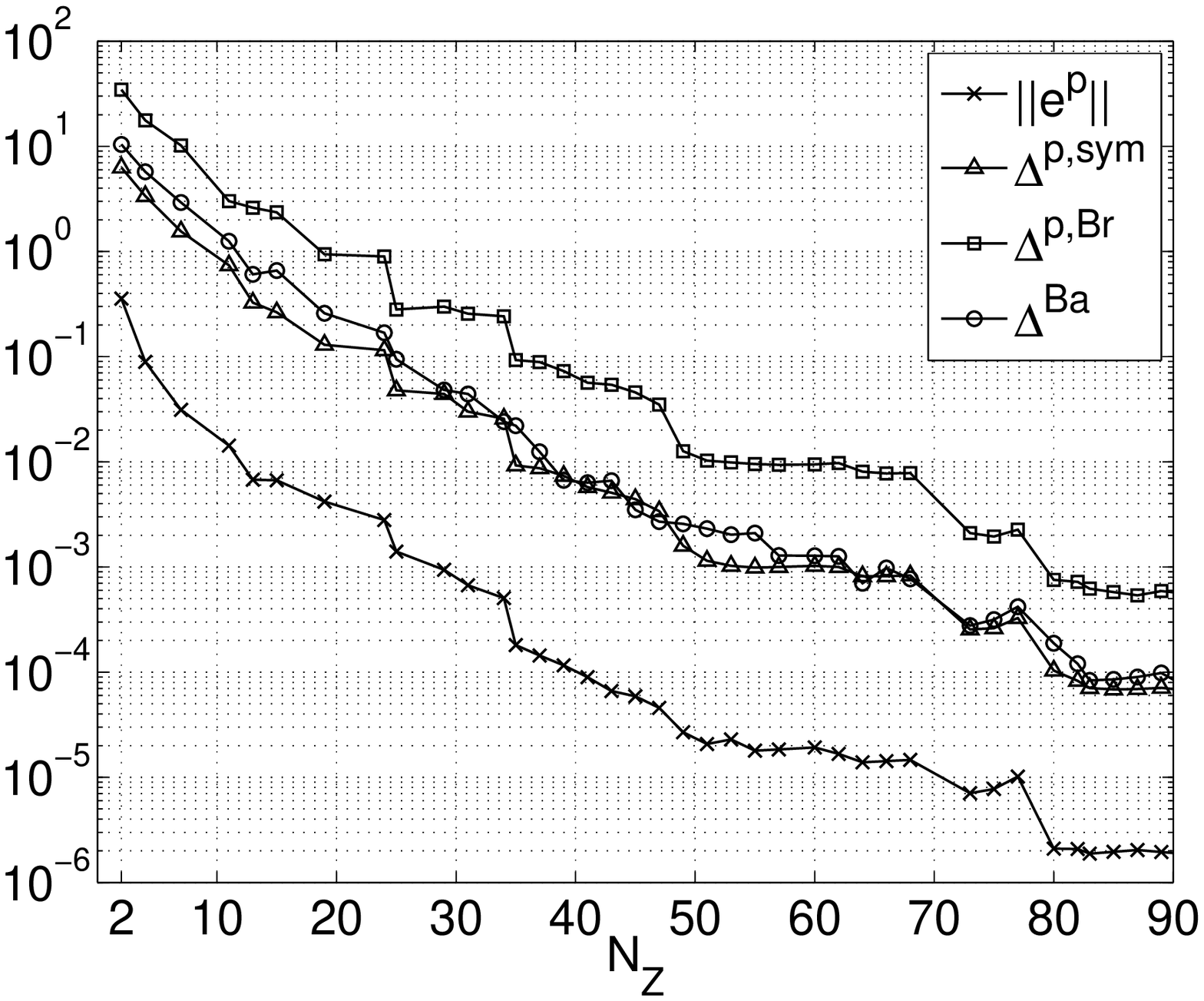,scale=0.3}\hfill
\centerline{(a) \hspace{35ex} (b) \hspace{35ex} (c)}
\caption{Maximum relative error $\|e^p_N(\mu)\|_Y/\|p(\mu)\|_Y$ and maximum relative error bounds $\Delta^{p,\rm sym}_N(\mu)/\|p(\mu)\|_Y$, $\Delta^{p}_N(\mu)/\|p(\mu)\|_Y$, and $\Delta^\Ba_N(\mu)/\|p(\mu)\|_Y$ shown as functions of $N_Z$ for (a)~Algorithm~1, (b)~Algorithm~2, and (c)~Algorithm~3 (see \S \ref{ss:construction_RB_spaces}); the maximum is taken over 25 parameter values; the computation of the error bounds is based on the exact constants involved.}
\label{fig:p_error}
\end{figure}

\begin{figure}[htp]
\epsfig{file=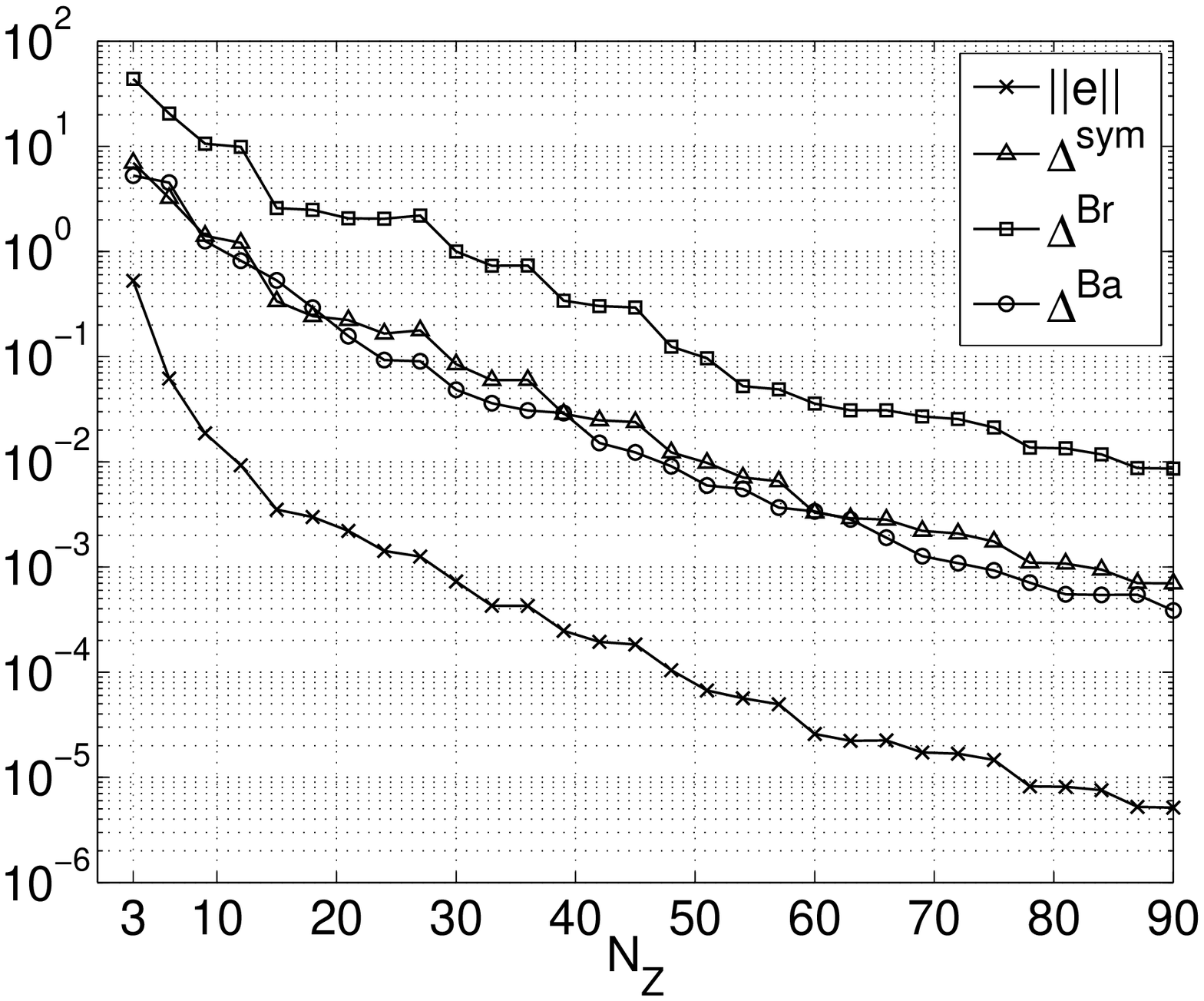,scale=0.3}\hfill
\epsfig{file=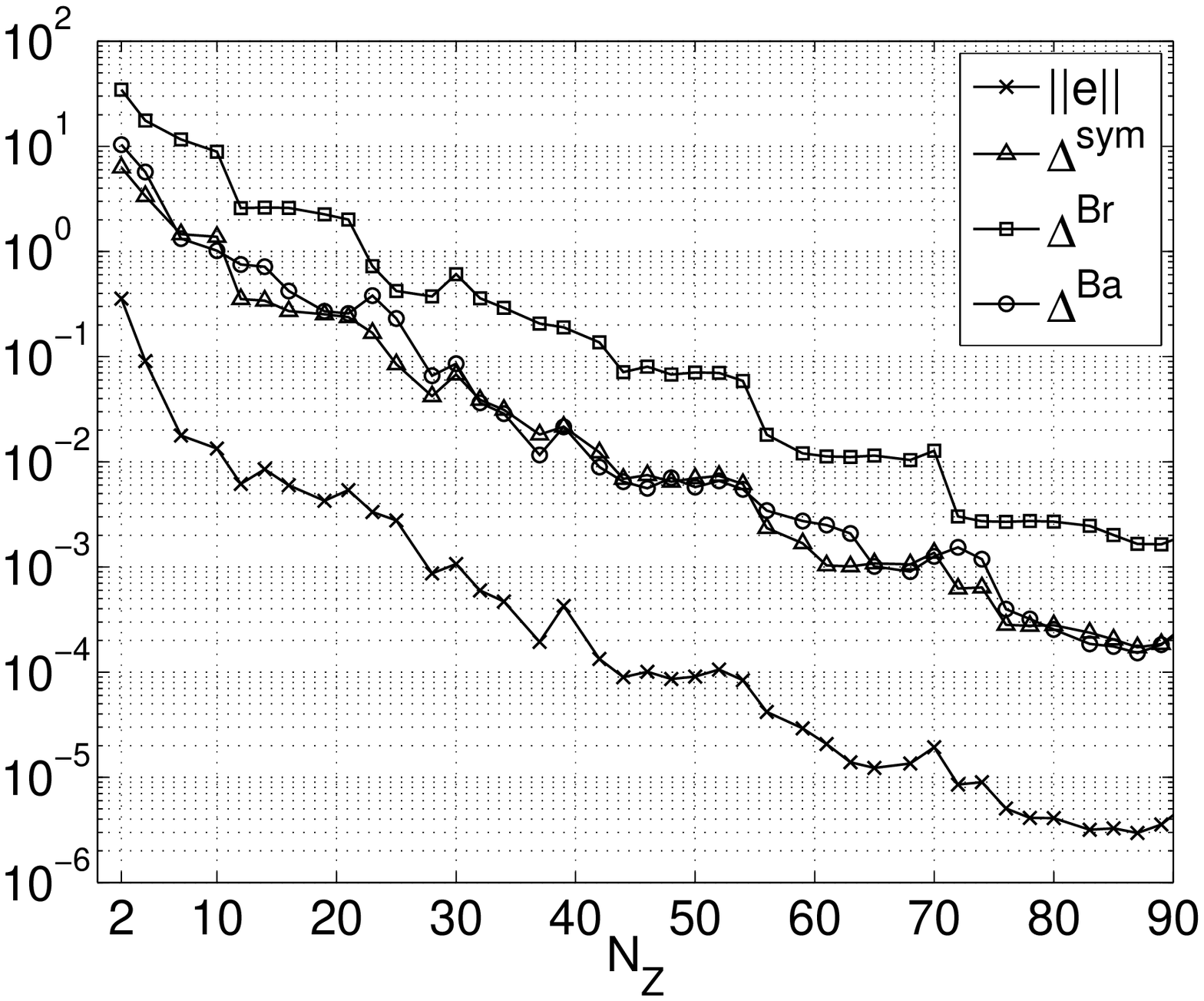,scale=0.3}\hfill
\epsfig{file=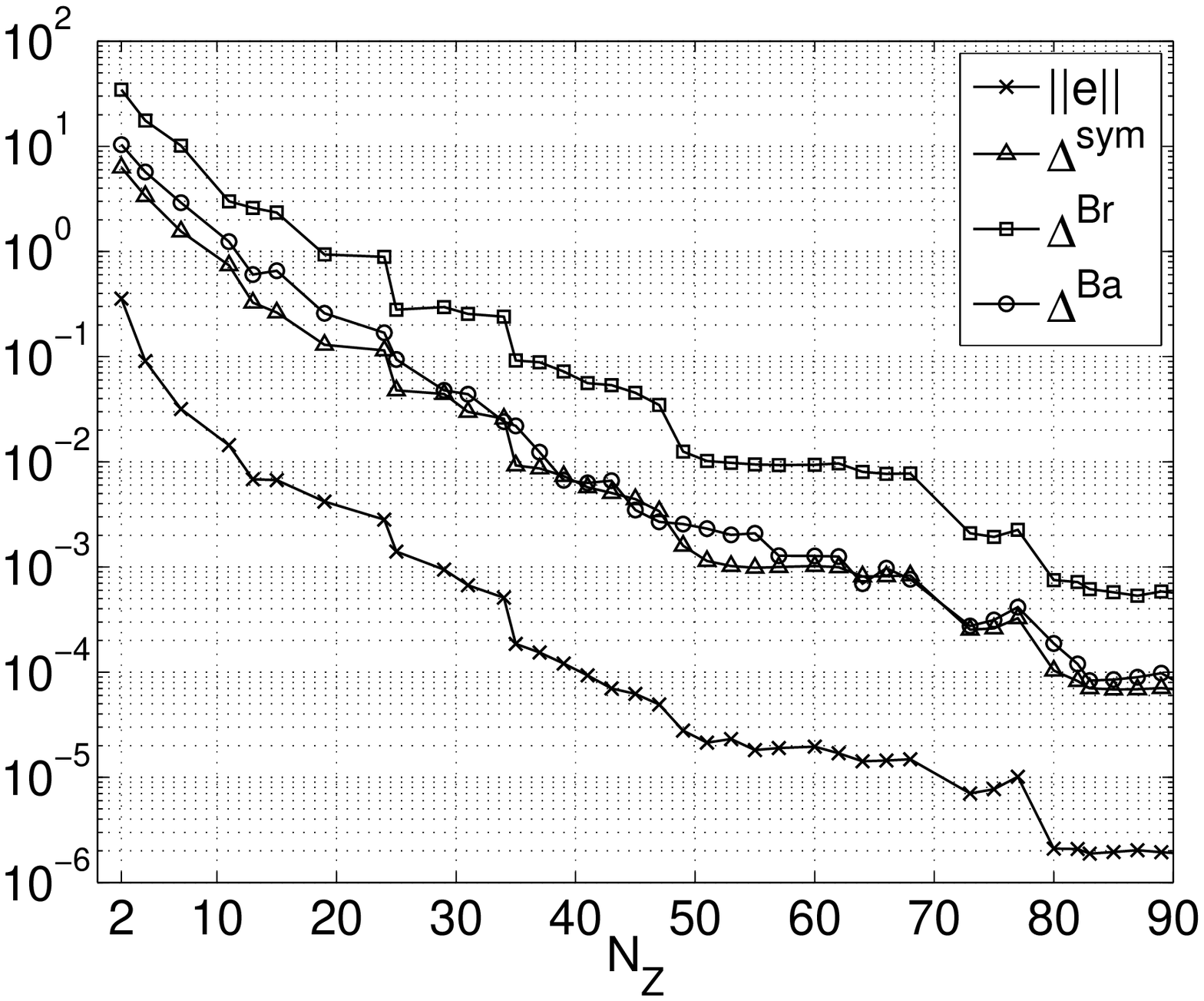,scale=0.3}\hfill
\centerline{(a) \hspace{35ex} (b) \hspace{35ex} (c)}
\caption{Maximum relative error $\|e_N(\mu)\|_Z/\|(u(\mu),p(\mu))\|_Z$ and maximum relative error bounds $\Delta^{\rm sym}_N(\mu)/\|(u(\mu),p(\mu))\|_Z$, $\Delta^\Br_N(\mu)/\|(u(\mu),p(\mu))\|_Z$, and $\Delta^\Ba_N(\mu)/\|(u(\mu),p(\mu))\|_Z$ shown as functions of $N_Z$ for (a)~Algorithm~1, (b)~Algorithm~2, and (c)~Algorithm~3 (see \S \ref{ss:construction_RB_spaces}); the maximum is taken over 25 parameter values; the computation of the error bounds is based on the exact constants involved.}
\label{fig:total_error}
\end{figure}

\begin{table}[htp]
\centering \footnotesize
\begin{tabular}{@{}c|c|c|c|c|c|c|c|c|c|c|c@{}}
\multicolumn{12}{l}{(a) Algorithm~1}\\[1ex]
\hline
& & & & & & & & & & &\\[-2.4ex]
$N$ & $N_Z$ & $\tilde\eta^{u,\rm sym}_N$ & $\eta^{u,\rm sym}_N$ &  $\eta^{u,\Br}_N$ & $\eta^{u,\Ba}_N$ & $\eta^{p,\rm sym}_N$ & $\eta^{p,\Br}_N$ & $\eta^{p,\Ba}_N$ & $\eta^{\rm sym}_N$ & $\eta^\Br_N$ & $\eta^\Ba_N$ \\[0.5ex]
\hline
$5$ & $15$ & $13.21$ & $20.68$ & $58.04$ & $353.6$ & $202.1$ & $1534$ & $296.4$ & $148.0$ & $1238$ & $203.6$\\
$10$ & $30$ & $12.34$ & $19.94$ & $63.93$ & $277.5$ & $292.5$ & $2268$ & $326.5$ & $168.6$ & $1559$ & $186.7$\\
$15$ & $45$ & $13.22$ & $19.57$ & $57.42$ & $309.9$ & $331.5$ & $3223$ & $417.2$ & $201.2$ & $1598$ & $228.6$\\
$20$ & $60$ & $13.29$ & $19.35$ & $58.40$ & $255.0$ & $355.9$ & $3240$ & $415.7$ & $208.1$ & $1580$ & $207.1$\\
$25$ & $75$ & $11.97$ & $18.77$ & $57.92$ & $242.9$ & $405.8$ & $3513$ & $476.3$ & $198.4$ & $1560$ & $216.4$\\
$30$ & $90$ & $13.55$ & $19.93$ & $54.23$ & $282.3$ & $399.2$ & $3423$ & $453.6$ & $227.0$ & $1678$ & $232.6$\\
$35$ & $105$ & $13.62$ & $20.44$ & $53.93$ & $281.7$ & $432.1$ & $2502$ & $519.8$ & $235.2$ & $1535$ & $234.9$\\
$40$ & $120$ & $14.28$ & $21.53$ & $57.49$ & $295.4$ & $482.6$ & $3317$ & $495.2$ & $253.0$ & $1526$ & $253.1$\\
\hline
\end{tabular}\\[1ex]
\begin{tabular}{@{}c|c|c|c|c|c|c|c|c|c|c|c@{}}
\multicolumn{12}{l}{(b) Algorithm~2}\\[1ex]
\hline
& & & & & & & & & & &\\[-2.4ex]
$N$ & $N_Z$ & $\tilde\eta^{u,\rm sym}_N$ & $\eta^{u,\rm sym}_N$ &  $\eta^{u,\Br}_N$ & $\eta^{u,\Ba}_N$ & $\eta^{p,\rm sym}_N$ & $\eta^{p,\Br}_N$ & $\eta^{p,\Ba}_N$ & $\eta^{\rm sym}_N$ & $\eta^\Br_N$ & $\eta^\Ba_N$ \\[0.5ex]
\hline
$5$ & $12$ & $16.11$ & $23.89$ & $55.18$ & $715.1$ & $98.47$ & $756.3$ & $142.5$ & $88.04$ & $673.4$ & $139.7$\\
$10$ & $23$ & $22.45$ & $32.53$ & $62.65$ & $1894$ & $121.8$ & $799.9$ &  $234.8$ & $99.72$ & $700.7$ & $233.0$\\
$15$ & $34$ & $13.57$ & $25.31$ & $68.15$ & $425.1$ & $152.8$ & $766.6$ & $253.9$ & $130.7$ & $660.6$ & $216.6$\\
$20$ & $46$ & $13.17$ & $20.31$ & $61.41$ & $505.3$ & $224.0$ & $1268$ & $307.0$ & $172.7$ & $1081$ & $236.1$\\
$25$ & $56$ & $16.44$ & $27.51$ & $63.04$ & $1013$ & $97.71$ & $550.1$ & $258.9$ & $94.87$ & $523.1$ & $250.8$\\
$30$ & $68$ & $14.63$ & $23.48$ & $61.14$ & $583.5$ & $237.9$ & $1287$ & $294.2$ & $189.6$ & $1023$ & $233.8$\\
$35$ & $78$ & $14.80$ & $22.35$ & $59.63$ & $707.3$ & $126.3$ & $843.6$ & $269.0$ & $118.3$ & $778.2$ & $251.5$\\
$40$ & $89$ & $14.13$ & $23.70$ & $61.52$ & $613.1$ & $162.4$ & $800.6$ & $285.4$ & $139.3$ & $685.6$ & $244.3$\\
\hline
\end{tabular}\\[1ex]
\begin{tabular}{@{}c|c|c|c|c|c|c|c|c|c|c|c@{}}
\multicolumn{12}{l}{(c) Algorithm~3}\\[1ex]
\hline
& & & & & & & & & & &\\[-2.4ex]
$N$ & $N_Z$ & $\tilde\eta^{u,\rm sym}_N$ & $\eta^{u,\rm sym}_N$ &  $\eta^{u,\Br}_N$ & $\eta^{u,\Ba}_N$ & $\eta^{p,\rm sym}_N$ & $\eta^{p,\Br}_N$ & $\eta^{p,\Ba}_N$ & $\eta^{\rm sym}_N$ & $\eta^\Br_N$ & $\eta^\Ba_N$ \\[0.5ex]
\hline
$5$ & $13$ & $18.36$ & $32.56$ & $68.23$ & $1146$ & $85.38$ & $793.3$ & $162.2$ & $79.09$ & $755.4$ & $157.9$\\
$10$ & $29$ & $22.08$ & $33.22$ & $72.73$ & $717.1$ & $183.0$ & $996.8$ & $219.0$ & $152.9$ & $830.9$ & $182.5$\\
$15$ & $39$ & $15.56$ & $26.58$ & $63.64$ & $446.0$ & $137.4$ & $721.9$ & $195.0$ & $123.8$ & $649.1$ & $173.1$\\
$20$ & $49$ & $15.13$ & $24.88$ & $69.96$ & $699.3$ & $118.1$ & $807.8$ & $174.9$ & $104.2$ & $753.7$ & $162.0$\\
$25$ & $60$ & $19.45$ & $33.86$ & $71.13$ & $1034$ & $130.3$ & $577.4$ & $206.5$ & $115.1$ & $509.1$ & $195.4$\\
$30$ & $73$ & $20.00$ & $40.18$ & $97.51$ & $976.6$ & $102.5$ & $582.7$ & $187.7$ & $94.86$ & $538.1$ & $179.6$\\
$35$ & $83$ & $16.90$ & $28.72$ & $71.93$ & $821.5$ & $136.3$ & $692.0$ & $217.4$ & $124.6$ & $631.3$ & $198.3$\\
$40$ & $93$ & $17.87$ & $33.31$ & $ 82.45$ & $973.4$ & $ 196.7$ & $1012$ & $301.9$ & $165.4$ & $849.5$ & $288.4$\\
\hline
\end{tabular}\\[2ex]
\caption{Maximum effectivities $\tilde\eta^{u,\rm sym}_N(\mu)=\tilde\Delta^{u,\rm sym}_N(\mu)/\|e^u_N(\mu)\|_{X,\mu}$ (see (\ref{eq:br_errbnd_u_sym_energy_rig})), $\eta^{u,\rm sym}_N(\mu)=\Delta^{u,\rm sym}_N(\mu)/\|e^u_N(\mu)\|_X$, $\eta^{p,\rm sym}_N(\mu)\equiv \Delta^{p,\rm sym}_N(\mu)/\|e^p_N(\mu)\|_Y$, $\eta^{\rm sym}_N(\mu)= \Delta^{\rm sym}_N(\mu)/\|e_N(\mu)\|_Z$ (see (\ref{eq:br_u_p_rig_sym}) and (\ref{eq:br_effectivity_sym})), and $\eta^{u,\Br}_N(\mu)\equiv \Delta^{u}_N(\mu)/\|e^u_N(\mu)\|_X$, $\eta^{p,\Br}_N(\mu)\equiv \Delta^{p}_N(\mu)/\|e^p_N(\mu)\|_Y$, $\eta^{\Br}_N(\mu)= \Delta^{\Br}_N(\mu)/\|e_N(\mu)\|_Z$, $\eta^{u,\Ba}_N(\mu)\equiv\Delta^\Ba_N(\mu)/\|e^u_N(\mu)\|_X$, $\eta^{p,\Ba}_N(\mu)\equiv \Delta^\Ba_N(\mu)/\|e^p_N(\mu)\|_Y$, $\eta^{\Ba}_N(\mu)= \Delta^\Ba_N(\mu)/\|e_N(\mu)\|_Z$ (see \cite{Gerner:2011fk}) for (a)~Algorithm~1, (b)~Algorithm~2, and (c)~Algorithm~3 (see \S \ref{ss:construction_RB_spaces}); the maximum is taken over 25 parameter values; the computation of the error bounds is based on the exact constants involved.}
\label{tbl:effectivities}
\end{table}

We now apply the RB methodology developed in \S \ref{s:RB_method} to the Stokes model problem described in \cite{Gerner:2012ag,Gerner:2011fk}. Numerical results are attained using the open source software \texttt{rbOOmit}~\cite{Knezevic:2011fk}, an implementation of the RB framework within the C++ parallel finite element library \texttt{libMesh}~\cite{libMeshPaper}. 

Motivated by applications in the field of microfluidics, we consider a Stokes flow in a two-dimensional microchannel with a parametrized, rectangular obstacle (see \cite{Gerner:2012ag}). As our truth discretization in \S \ref{ss:truth}, we then choose a mixed finite element method using the standard conforming $\mathbb{P}_2$-$\mathbb{P}_1$ Taylor--Hood approximation spaces \cite{Taylor:1973fk}; the truth system (\ref{eq:truth}) here exhibits a dimension of $\cN=$ 72,076. 

Since our Stokes model problem is clearly symmetric, we may now compare the RB {\em a posteriori} error bounds $\Delta^{u,\rm sym}_N(\mu)$, $\Delta^{p,\rm sym}_N(\mu)$, and $\Delta^{\rm sym}_N(\mu)$ (see (\ref{eq:br_errbnd_u_sym}), (\ref{eq:br_errbnd_p_sym}), and (\ref{eq:br_errbnd_sym})) specific to the symmetric case with the more general bounds $\Delta^{u}_N(\mu)$, $\Delta^{p}_N(\mu)$, and $\Delta^\Br_N(\mu)$ introduced in \cite[Corollary~2.2]{Gerner:2011fk}; in addition, we shall again consider the RB error estimate $\Delta^\Ba_N(\mu)$ based on Babu\v{s}ka's theory for noncoercive problems as defined in \cite[Corollary~2.1]{Gerner:2011fk}. For this purpose, we build the approximation spaces $X_N$ and $Y_N$ by using either Algorithm~1, Algorithm~2, or Algorithm~3 (see \S \ref{ss:construction_RB_spaces}). All sampling procedures are based on an exhaustive random sample $\Sigma\subset\cD$ of size $|\Sigma|=$ 4,900 and the relative error bound $\Delta_N(\mu)\equiv \tilde\Delta^{u,\rm sym}_N(\mu)/\|u_N(\mu)\|_X$ (see (\ref{eq:br_errbnd_u_sym_energy})) for which (\ref{eq:eff_sym}) suggests the lowest effectivities; in Algorithm~2 and Algorithm~3, we further set $\delta^\beta_{\rm tol}=0.1$.

Figure~\ref{fig:u_error} and Figure~\ref{fig:p_error} show the maximum errors in the RB approximations for the primal variable and the Lagrange multiplier, respectively, together with the respective error bounds $\Delta^{u,\rm sym}_N(\mu)$, $\Delta^{u}_N(\mu)$, $\Delta^{p,\rm sym}_N(\mu)$, $\Delta^{p}_N(\mu)$, and $\Delta^\Ba_N(\mu)$; Figure~\ref{fig:total_error} shows the maximum total error in the RB approximation and associated error bounds $\Delta^{\rm sym}_N(\mu)$, $\Delta^\Br_N(\mu)$, and $\Delta^\Ba_N(\mu)$. The maximum is computed over a sample of 25 parameter values. For this sample, to analyze only the effects of the {\em a posteriori} error bound formulations and eliminate contributions of the SCM (see \S \ref{ss:offline_online}), we use the exact constants rather than the lower/upper bounds (\ref{eq:alpha_gamma_LB}). 
Effectivities associated with the error bounds are given in Table~\ref{tbl:effectivities}. As before (see \cite{Gerner:2011fk}), maximum effectivities associated with $\Delta^{u}_N(\mu)$ essentially range from $50$ to $80$ (see Fig.~\ref{fig:u_error} and Table~\ref{tbl:effectivities}). Now, exploiting the symmetry of the problem, $\Delta^{u,\rm sym}_N(\mu)$ provides a sharper bound not only in theory (see Proposition~\ref{prpstn:br_errbnds_sym}) but also in practice: Here, effectivities essentially vary between $20$ and $40$. Furthermore, it is not surprising (see (\ref{eq:eff_sym})) to observe that the best results are achieved by the bound $\tilde\Delta^{u,\rm sym}_N(\mu)$ that overestimates the error in the RB approximations for the primal variable measured in the energy norm (see (\ref{eq:br_errbnd_u_sym_energy_rig})) only by a factor of approximately $15$ (see Table~\ref{tbl:effectivities}).
For $\Delta^{p,\rm sym}_N(\mu)$, we still obtain rather large effectivities with average values of 362, 161, and 127 in case of Algorithm~1, Algorithm~2, and Algorithm~3, respectively; however, compared to $\Delta^{p}_N(\mu)$, the improvement is significant (see Fig.~\ref{fig:p_error} and Table~\ref{tbl:effectivities}). We emphasize that $\Delta^{u,\rm sym}_N(\mu)$ and $\Delta^{p,\rm sym}_N(\mu)$ do in neither case perform worse than $\Delta_N^\Ba(\mu)$ but generally provide bounds for the errors in the primal and Lagrange multiplier variables that are much more accurate (see Table~\ref{tbl:effectivities}).

We now discuss the Online computation times of the proposed methods. 
The SCM (see \S \ref{ss:offline_online}) enables the (Online-)efficient estimation of the constants $\alpha_a(\mu)$, $\gamma_a(\mu)$, and $\beta_\Br(\mu)$. We here apply the method with the configurations as specified in  \cite{Gerner:2011fk} and receive very accurate (Online-)efficient bounds $\alpha^\LB_a(\mu)$, $\gamma^\UB_a(\mu)$, and $\beta^\LB_\Br(\mu)$, providing {\em a posteriori} error bounds $\Delta^{u}_N(\mu)$, $\Delta^{u,\rm sym}_N(\mu)$, $\tilde\Delta^{u,\rm sym}_N(\mu)$ and $\Delta^{p}_N(\mu)$, $\Delta^{p,\rm sym}_N(\mu)$ that essentially coincide with their values based on the evaluation of the exact constants; associated effectivities thus remain the same as shown in Table~\ref{tbl:effectivities}. For comparison, once the $\mu$-independent parts in the affine expansion of the involved operators (see \S \ref{ss:offline_online}) have been formed, direct computation of the truth approximation $(u(\mu),p(\mu))$ (i.e., assembly and solution of (\ref{eq:truth})) takes on average 6.5~seconds on a 2.66 GHz Intel Core 2 Duo processor. The rigorous and (Online-)efficient error bounds $\Delta^{u,\rm sym}_N(\mu)$ and $\Delta^{p,\rm sym}_N(\mu)$ allow us to choose the RB system dimension $N_Z$ just large enough to obtain a desired accuracy. In case of Algorithm~1, we need $N_Z=51$ to achieve a prescribed accuracy of roughly $1\%$ or better in the RB approximations $u_N(\mu)$ (see Fig.~\ref{fig:u_error}(a)). Once the database has been loaded, the Online calculation of $(u_N(\mu),p_N(\mu))$ (i.e., assembly and solution of (\ref{eq:rb})) and $\Delta^{u,\rm sym}_N(\mu)$, $\Delta^{p,\rm sym}_N(\mu)$ for any new value of $\mu\in\cD$ takes on average $0.31$ and $20.99$ milliseconds, respectively, which is in total roughly $300$ times faster than direct computation of the truth approximation. 
We note that again (see \cite{Gerner:2011fk}), stabilizing adaptively pays off:
In the case of Algorithm~3, the same accuracy is achieved for $N_Z=35$ (see Fig.~\ref{fig:u_error}(c)); the Online calculation of $(u_N(\mu),p_N(\mu))$ and $\Delta^{u,\rm sym}_N(\mu)$, $\Delta^{p,\rm sym}_N(\mu)$ then takes on average $0.14$ and $13.29$ milliseconds, respectively, and is thus roughly $480$ times faster than direct computation of the truth approximation.
Detailed computation times, also for Algorithm~2, are given in Table~\ref{tbl:computation_times}. Clearly, due to the significant improvement of the new error bounds in terms of sharpness, the new methods guarantee a prescribed accuracy in the RB approximations at notable Online savings when compared to the methods presented in \cite{Gerner:2011fk} (see Table~\ref{tbl:computation_times} and \cite[Table~2.3]{Gerner:2011fk}).

\begin{table}[tp]
\centering \footnotesize
\begin{tabular}{@{}c|c|c|c|c|c@{}}
\multicolumn{6}{l}{(a) Accuracy of at least 1\% (resp., 0.1\%) for the RB approximations $u_N(\mu)$}\\[1ex]
\hline
& & & & &\\[-2.4ex]
Method & $N_Z$ &  $N$ & $(u_N(\mu),p_N(\mu))$ & $\Delta^{u,\rm sym}_N(\mu), \Delta^{p,\rm sym}_N(\mu)$ & Total\\
& & & & &\\[-2.4ex]
\hline
Algorithm~1 & 51 (78) & 17 (26) & 0.31 (0.79) & 20.99 (39.60) & 21.20 (40.38)\\
\hline
Algorithm~2 & 44 (61) & 19 (27) & 0.21 (0.42) & 16.25 (24.73) & 16.46 (25.15)\\
\hline
Algorithm~3 & 35 (55) & 13 (23) & 0.14 (0.32) & 13.29 (21.59) & 13.43 (21.92)\\
\hline
\end{tabular}\\[1ex]
\begin{tabular}{@{}c|c|c|c|c|c@{}}
\multicolumn{6}{l}{(b) Accuracy of at least 1\% (resp., 0.1\%) for the RB approximations $p_N(\mu)$}\\[1ex]
\hline
& & & & &\\[-2.4ex]
Method & $N_Z$ &  $N$ & $(u_N(\mu),p_N(\mu))$ & $\Delta^{u,\rm sym}_N(\mu), \Delta^{p,\rm sym}_N(\mu)$ & Total\\
& & & & &\\[-2.4ex]
\hline
Algorithm~1 & 54 (84) & 18 (28) & 0.35 (0.91) & 22.75 (44.70) & 23.10 (45.61)\\
\hline
Algorithm~2 & 44 (72) & 19 (32) & 0.21 (0.58) & 16.25 (31.54) & 16.46 (32.12)\\
\hline
Algorithm~3 & 37 (62) & 14 (26) & 0.16 (0.44) & 13.94 (25.77) & 14.09 (26.21)\\
\hline
\end{tabular}\\[2ex]
\caption{Average computation times in milliseconds for the Online evaluation of $(u_N(\mu),p_N(\mu))$ (assembly and solution of (\ref{eq:rb})) and the error bounds $\Delta^{u,\rm sym}_N(\mu)$, $\Delta^{p,\rm sym}_N(\mu)$ (see (\ref{eq:br_errbnd_u_sym}), (\ref{eq:br_errbnd_p_sym})); times are measured using either Algorithm~1, Algorithm~2, or Algorithm~3 (see \S \ref{ss:construction_RB_spaces}), with a prescribed accuracy of 1\% (resp., 0.1\%) for the RB approximations (a) $u_N(\mu)$ and (b) $p_N(\mu)$.}
\label{tbl:computation_times}
\end{table}

\section{Conclusion} \label{s:conclusion}

In this paper, we improve RB error bounds introduced in \cite{Gerner:2011fk} for the special case of a symmetric problem. Numerical results provide a direct comparison with former approaches and demonstrate the superiority of the developed bounds with respect to sharpness; in particular, the upper bounds provided for the errors in the RB approximations for the primal variable exhibit effectivities that are comparatively low. 
As a direct consequence, the proposed methods provide accurate RB approximations at much less computational cost.

The analysis presents possible techniques but clearly does not claim to be exhaustive. For example, as current RB techniques allow their exact computation, {\em a posteriori} error bounds are here exclusively formulated in terms of the residual dual norms (\ref{eq:residual_dual_norms}). However, the symmetric case allows us to consider the dual norm of the residual (\ref{eq:residual_1}) with respect to the energy norm $\|\cdot\|_{X,\mu} = \sqrt{a(\cdot,\cdot;\mu)}$,
\begin{equation*}
\|r^1_N(\cdot;\mu)\|_{X',\mu}\equiv \sup_{v\in X}\frac{r^1_N(v;\mu)}{\|v\|_{X,\mu}} \quad\forall\;\mu\in\cD.
\end{equation*}
Using the same techniques as presented in the proof of Proposition~\ref{prpstn:br_errbnds_sym}, the errors $e^u_N(\mu)$ and $e^p_N(\mu)$ in the RB approximations for the primal and Lagrange multiplier variables may then be bounded in terms of $\|r^1_N(\cdot;\mu)\|_{X',\mu}$, $\|r^2_N(\cdot;\mu)\|_{Y'}$, and $\tilde\beta(\mu)$ (see (\ref{eq:LBB_tilde})): For any $\mu\in\cD$ and $N\in\bNmax$, we can derive that
\begin{align}\label{eq:spp_energy_errbnd_u}
\|e^u_N(\mu)\|_{X,\mu} &\leq \|r^1_N(\cdot;\mu)\|_{X',\mu} + \frac{1}{\tilde\beta(\mu)}\|r^2_N(\cdot;\mu)\|_{Y'},\\ \label{eq:spp_energy_errbnd_p}
\|e^p_N(\mu)\|_Y &\leq \frac{2}{\tilde\beta(\mu)}\|r^1_N(\cdot;\mu)\|_{X',\mu} + \frac{1}{(\tilde \beta(\mu))^2}\|r^2_N(\cdot;\mu)\|_{Y'}.
\end{align}
For these to provide useful error bounds in the RB context, $\|r^2_N(\cdot;\mu)\|_{X',\mu}$ and $\tilde\beta(\mu)$ need to be estimated Online-efficiently. Through (\ref{eq:alpha_a}) and (\ref{eq:tilde_beta_equiv}), this may be done in terms of $\|r^1_N(\cdot;\mu)\|_{X'}$, $\alpha_a(\mu)$, $\gamma_a(\mu)$, and $\beta_\Br(\mu)$; in this case, (\ref{eq:spp_energy_errbnd_u}) leads to (\ref{eq:br_errbnd_u_sym}) as well and (\ref{eq:spp_energy_errbnd_p}) results in an upper bound worse than (\ref{eq:br_errbnd_p_sym}). 

\begin{acknowledgement}
The authors are most grateful to Prof.~Arnold Reusken of RWTH Aachen University for numerous very helpful comments and suggestions. The authors would also like to thank Dr.~David J. Knezevic of Harvard University for his reliable support on {\tt rbOOmit}~\cite{Knezevic:2011fk}.
\end{acknowledgement}

\bibliographystyle{siam}
\bibliography{bibtex/General,bibtex/Reduced_Basis}

\end{document}